\documentclass{IEEEtran}

\usepackage{enumitem,kantlipsum}
\usepackage{graphicx}
\usepackage{epstopdf}
\usepackage{amssymb}
\usepackage{amsmath}
\usepackage{subfigure}
\usepackage{epsfig}
\usepackage{color}
\usepackage{enumitem}
\usepackage{float}
\usepackage{soul}
\usepackage{wrapfig}
\usepackage{arydshln}
\usepackage{tablefootnote}

\usepackage{amsthm}

\def\0{{\bf 0}}
\def\1{{\bf 1}}

\def\beq{\begin{equation*}}
    \def\eeq{\end{equation*}}
\def\bql{\begin{equation}}
    \def\eql{\end{equation}}
\def\bqn{\begin{eqnarray*}}
    \def\eqn{\end{eqnarray*}}
\def\bnl{\begin{eqnarray}}
    \def\enl{\end{eqnarray}}
\def\bma{\begin{bmatrix}}
    \def\ema{\end{bmatrix}}
\def\bmx{\begin{matrix}}
    \def\emx{\end{matrix}}
\def\ben{\begin{enumerate}}
    \def\een{\end{enumerate}}
\def\bit{\begin{itemize}}
    \def\eit{\end{itemize}}
\def\bei{\begin{itemize}}
    \def\eei{\end{itemize}}
\def\bet{\begin{tabular}}
    \def\eet{\end{tabular}}

\newcommand{\ba}{\mathbf{a}}

\newcommand{\R}{\mathbb{R}}

\newcommand{\be}{\mathbf{e}}

\newcommand{\g}{\mathbf{g}}

\newcommand{\bu}{\mathbf{u}}

\newcommand{\bv}{\mathbf{v}}

\usepackage[dvipsnames]{xcolor}

\def\R{\mathbb{R}}

\def\1{{\bf1}}
\def\la{\langle}
\def\ra{\rangle}
\def\b{{\beta}}
\def\a{\alpha}
\def\g{\gamma}

\def\bit{\begin{itemize}}
\def\eit{\end{itemize}}

\def\be{\begin{equation}}
\def\ee{\end{equation}}
\def\ba{\begin{eqnarray}}
\def\ea{\end{eqnarray}}

\def\bes{\begin{equation*}}
\def\ees{\end{equation*}}
\def\bas{\begin{eqnarray*}}
\def\eas{\end{eqnarray*}}

\usepackage{url}
\usepackage{verbatim}

\newtheorem{Remark 1}{Remark}
\newtheorem{Remark 2}[Remark 1]{Remark}
\newtheorem{Remark 3}[Remark 1]{Remark}
\newtheorem{Remark 4}[Remark 1]{Remark}
\newtheorem{Remark 5}[Remark 1]{Remark}
\newtheorem{Remark 6}[Remark 1]{Remark}
\newtheorem{Remark 7}[Remark 1]{Remark}

\newtheorem{Lemma 1}{Lemma}
\newtheorem{Lemma 2}[Lemma 1]{Lemma}
\newtheorem{Lemma 3}[Lemma 1]{Lemma}
\newtheorem{Lemma 4}[Lemma 1]{Lemma}
\newtheorem{Lemma 5}[Lemma 1]{Lemma}
\newtheorem{Lemma 6}[Lemma 1]{Lemma}
\newtheorem{Lemma 7}[Lemma 1]{Lemma}

\newtheorem{Assumption 1}{Assumption}
\newtheorem{Assumption 2}[Assumption 1]{Assumption}
\newtheorem{Assumption 3}[Assumption 1]{Assumption}
\newtheorem{Assumption 4}[Assumption 1]{Assumption}
\newtheorem{Definition 1}{Definition}
\newtheorem{Theorem 1}{Theorem}
\newtheorem{Theorem 2}[Theorem 1]{Theorem}
\newtheorem{Theorem 3}[Theorem 1]{Theorem}
\newtheorem{Theorem 4}[Theorem 1]{Theorem}
\newtheorem{Theorem 5}[Theorem 1]{Theorem}
\newtheorem{Theorem 6}[Theorem 1]{Theorem}
\newtheorem{Theorem 7}[Theorem 1]{Theorem}
\newtheorem{Theorem 8}[Theorem 1]{Theorem}
\newtheorem{Theorem 9}[Theorem 1]{Theorem}
\newtheorem{Theorem 10}[Theorem 1]{Theorem}

\title{\LARGE \bf
 Tailoring Gradient Methods for Differentially-Private Distributed Optimization
 }

\author{Yongqiang Wang, Angelia Nedi\'c
\thanks{ The work was supported in part by the National Science Foundation under Grants   ECCS-1912702, CCF-2106293, and CCF-2106336.}
\thanks{Yongqiang Wang is with the Department of Electrical and Computer Engineering, Clemson University, Clemson, SC 29634, USA
{\tt\small{yongqiw}@clemson.edu}
}%
\thanks{Angelia Nedi\'c is with the School of Electrical, Computer and Energy Engineering, Arizona State University, Tempe, AZ 85281, USA {\tt\small angelia.nedich@asu.edu}}
}

\begin{document}

\maketitle
\thispagestyle{empty}
\pagestyle{empty}

\begin{abstract}
Decentralized optimization is gaining increased traction due to its widespread applications in large-scale machine learning and multi-agent systems. The same mechanism that enables its success, i.e., information sharing among participating agents, however, also leads to the
disclosure of individual agents' private information, which is unacceptable when  sensitive data are involved.
As differential privacy is becoming a de facto standard for privacy preservation, recently  results have emerged integrating differential privacy with distributed optimization.
However, directly incorporating differential privacy design in existing distributed optimization approaches   significantly compromises  optimization accuracy. In this paper, we propose to redesign and tailor  gradient methods for differentially-private distributed optimization, and propose two differential-privacy oriented gradient methods that can ensure both {rigorous $\epsilon$-differential privacy} and optimality.
%
The first algorithm is based on static-consensus based gradient methods, and the second algorithm is based on dynamic-consensus (gradient-tracking) based distributed optimization methods and, hence,  is applicable to general directed interaction  graph  topologies.
{Both  algorithms can  simultaneously ensure   almost sure convergence  to an optimal solution and a finite privacy budget, even when the number of iterations goes to infinity. To our knowledge, this is the first time that  both goals are achieved simultaneously.}   Numerical {simulations using a distributed estimation problem and experimental results on a benchmark dataset} confirm  the effectiveness of the proposed approaches.
\end{abstract}

\section{Introduction}

The problem of  optimizing a global objective function through the cooperation of multiple agents has gained increased attention in recent years. This is driven by its wide applicability to many engineering and scientific domains, ranging from cooperative control
\cite{yang2019survey}, distributed sensing
\cite{bazerque2009distributed},  sensor networks
\cite{zhang2017distributed}, to large-scale machine learning
\cite{tsianos2012consensus}. In many of these applications, each agent only has access to a local objective function and can only communicate with its local neighbors. The agents cooperate to minimize the summation of all individual agents' local objective functions. Such a distributed optimization problem can be formulated in the following general form:
\begin{equation}\label{eq:optimization_formulation1}
\min\limits_{\theta\in\mathbb{R}^d} F(\theta)\triangleq
\frac{1}{m}\sum_{i=1}^m f_i(\theta)
\end{equation}
where $m$ is the number of agents, $\theta\in\mathbb{R}^d$ is
 a decision variable common to all agents, while
$f_i:\mathbb{R}^d\rightarrow\mathbb{R}$ is a local objective
function private to agent $i$.

Plenty of approaches have been reported to solve the above
distributed optimization problem since the seminal work of
\cite{tsitsiklis1984problems}, with some of the commonly used approaches
including gradient methods (e.g.,
\cite{nedic2009distributed,srivastava2011distributed,shi2015extra,xu2017convergence,qu2017harnessing,xin2018linear}),
distributed alternating direction method of multipliers (e.g.,
\cite{shi2014linear,zhang2019admm}),  and distributed Newton methods
(e.g., \cite{wei2013distributed}).
Among   these approaches, gradient-based approaches are gaining increased traction due to
their efficiency in both computation complexity and storage requirement, which is particularly appealing for   agents with limited computational or storage capabilities.
In general, existing gradient based distributed optimization algorithms can be divided into two categories. The first category  combines gradient-descent operations and average-consensus mechanisms
(referred to as  static-consensus hereafter) by directly concatenating gradient-descent with a consensus operation of individual agents' optimization variables. Typical examples include~\cite{nedic2009distributed,yuan2016convergence}. Such approaches are simple and efficient in computation since they only require an agent to share one variable in each iteration. However,
 these  approaches  are only applicable in
balanced graphs (the sum of each agent's in-neighbor coupling weights equal to the sum of its out-neighbor coupling weights). The second category circumvents the balanced-graph restriction by exploiting consensus mechanisms able to track time-varying signals (so-called dynamic consensus, applicable to general directed graphs) to track  the global gradient (see, e.g., \cite{xu2017convergence,qu2017harnessing,xin2018linear,di2016next,pu2020push}). It can ensure convergence to an optimal solution under constant  stepsizes and, hence, can achieve faster convergence. However, such approaches need every agent to maintain and share an additional gradient-tracking variable besides the optimization variable, which doubles the communication overhead.

Despite the enormous success of  gradient based distributed optimization algorithms, they all explicitly share optimization variables and/or gradient estimates  in every iteration, which becomes a problem in applications involving sensitive data. For example, in
the rendezvous problem where a group of agents uses distributed
optimization to cooperatively find an optimal assembly point,
participating agents may want to keep their initial positions
private, which is particularly important in unfriendly
environments~\cite{zhang2019admm}.
In sensor network based localization, the positions of sensor agents should be kept private  in sensitive (hostile) environments as well~\cite{zhang2019admm,huang2015differentially}.
In fact, without an effective privacy
mechanism in place, the results
in~\cite{zhang2019admm,huang2015differentially,burbano2019inferring}
show  that a participating agent's sensitive information, such as
position, can be easily  inferred by an adversary or other
participating agents in distributed-optimization based rendezvous
and localization approaches.
Another example underscoring the importance
of privacy protection in distributed optimization is machine
learning where exchanged data may contain sensitive information such
as medical records or salary information~\cite{yan2012distributed}.
In fact, recent  results in~\cite{zhu2019deep} show  that without a
privacy mechanism in place, an adversary can use shared information to precisely recover the raw data used for training
(pixel-wise accurate for images and token-wise matching for texts).

To address the pressing need for privacy protection in distributed optimization, recently plenty of efforts have been reported to
counteract potential privacy breaches in distributed optimization.
One approach resorts to partially homomorphic encryption, which has been employed in both our own prior results~\cite{zhang2019admm,zhang2018enabling}, and
others~\cite{freris2016distributed,lu2018privacy}.
However, such approaches
   incur  heavy communication and computation overhead.
   Another approach employs  the structural
  properties of distributed optimization to inject temporally or spatially correlated uncertainties,  which can also provide privacy protection in distributed optimization. For example,
\cite{yan2012distributed,lou2017privacy} showed that privacy can be enabled by
  adding a {\it constant} uncertain parameter in the projection step or stepsizes. The authors of \cite{gade2018private} showed that network structure can be leveraged to
  construct spatially correlated ``structured" noise to cover
  information.
   However,
since the uncertainties injected by these approaches are correlated, their enabled privacy is restricted: projection based privacy depends on
  the size of the projection set -- a large projection set nullifies
  privacy protection whereas a small projection set offers strong
  privacy
  protection but requires {\it a priori} knowledge of the optimal solution; ``structured"
  noise based approaches require  each agent to have a certain number
  of neighbors that do not share information with the adversary.
Differential Privacy (DP)~\cite{dwork2014algorithmic} is becoming increasingly  popular in privacy protection. It employs uncorrelated noises, and hence can provide strong privacy protection for a participating agent, even when all its neighbors  are compromised.
  As DP  is achieving remarkable successes in various
applications~\cite{han2016differentially,hale2017cloud,wang2017differential,zhang2019recycled,he2020differential} and  becoming a de facto standard for privacy protection, some efforts
have also been reported incorporating DP-noise into distributed optimization.
  For example,
 approaches have been proposed  to obscure
shared information in distributed optimization by injecting DP-noise
to exchanged messages~\cite{huang2015differentially,cortes2016differential,xiong2020privacy,ding2021differentially},
or objective functions~\cite{nozari2016differentially}. However,
while obscuring information, directly incorporating persistent DP-noise into existing algorithms also unavoidably compromises the accuracy of optimization, leading to a fundamental trade-off between privacy and
accuracy. In fact, recently the investigation in~\cite{zhu2019deep} indicates that DP-based defense can achieve reasonable privacy protection {\it ``only when the noise variance is
large enough to degrade accuracy~\cite{zhu2019deep}."}

In this paper, we propose to  tailor  gradient methods for differentially-private distributed  optimization. More specifically,  motivated by the observation that persistent DP-noise has to be repeatedly injected in every iteration of gradient based methods to ensure a strong privacy protection, which results in significant reduction in optimization accuracy, we propose to  gradually weaken coupling strength in distributed optimization to attenuate DP-noise that is added to every shared message.
We  judiciously design  the weakening factor sequences such that the consensus and convergence
to an  optimal solution are ensured even in the presence of persistent DP-noise.

The main contributions are as follows:
1)~We propose two
 gradient-based  methods for differentially private distributed optimization.
 The first one is based on static-consensus combined with a gradient method, which needs every agent to store and share one variable in each iteration. The second one
is based on dynamic-consensus (gradient-tracking) combined with an approximate gradient method,  which needs every agent to store and share two variables, but it is applicable to general directed graphs;
2)~We rigorously prove that both algorithms can ensure almost sure convergence of all agents to the optimal solution even in the presence of persistent DP-noise, which, to our knowledge, has not been  achieved before;
3)~{ We   prove that both algorithms can ensure rigorous $\epsilon$-DP for participating agents' objective functions, even when all communications are observable to adversaries. More interestingly, { both algorithms can  ensure a finite privacy budget even when the number of iterations goes to infinity. To our knowledge, this is the first time that  almost sure convergence to the optimal solution and rigorous $\epsilon$-DP (with a guaranteed finite privacy budget even when the number of iterations tends to infinity) are achieved simultaneously in distributed optimization}; 4) Even without taking privacy into consideration, the two proposed algorithms and theoretical derivations are of interest themselves.
 We propose a new vector-valued martingale convergence theorem (Lemma \ref{th-dsystem}) as well as its adaptations to  distributed optimization problems (Lemmas \ref{th-main_decreasing}, \ref{Theo:convergence_to_zero}, and \ref{Theorem:general_gradient_tracking}), which enables us to analyze the consensus-error evolution and optimality-gap evolution under DP-noise simultaneously. }

The organization of the paper is as follows.
Sec.~\ref{sec-problem} gives the problem formulation
and some results for a later use.
Sec.~\ref{sec-static} presents  a static-consensus based  gradient method for differentially-private distributed optimization and  establishes the almost sure  convergence of all agents' iterates to an optimal solution {as well as  $\epsilon$-DP guarantees}.
Sec.~\ref{sec-dynamic} presents
 a dynamic-consensus based gradient  method for differentially-private distributed optimization and establishes the almost sure  convergence   to an optimal solution {as well as   $\epsilon$-DP guarantees}.   Sec.~\ref{sec-numerics} presents   both { numerical simulations and experimental results on a benchmark dataset MNIST}. Finally, Sec.~\ref{sec-concl}  concludes the paper.

{\bf Notations:}  We use $\mathbb{R}^d$ to denote the Euclidean space of
dimension $d$. We write $I_d$ for the identity matrix of dimension $d$,
and ${\bf 1}_d$ for  the $d$-dimensional  column vector will all
entries equal to 1; in both cases  we suppress the dimension when
clear from the context. 
For a
vector $x$, $x_i$ denotes its $i$th element.
  We use $\langle\cdot,\cdot\rangle$ to denote the inner product. 
We write $\|A\|$ for the matrix norm induced by the vector norm $\|\cdot\|$, unless stated otherwise.
We let $A^T$ denote the transpose of a matrix $A$.
We also use other vector/matrix norms defined under a certain transformation determined by a matrix $W$, which will be represented as $\|\cdot\|_W$.
A matrix is column-stochastic when its entries are nonnegative and
elements in every column add up to one. A square matrix $A$ is said
to be doubly-stochastic when both $A$ and $A^T$ are
column-stochastic. For two vectors
$u$ and $v$ with the same dimension, we use $u \leq v$ to represent
the relationship that every element of the vector $u-v$ is nonpositive. Often, we abbreviate {\it almost surely} by {\it a.s}.
\def\as{{\it a.s.\ }}

\section{Problem Formulation and Preliminaries}\label{sec-problem}

\subsection{On distributed optimization}
We consider  a network of $m$ agents, interacting on a general directed graph. We describe a directed graph using an ordered pair $\mathcal{G}=([m],\mathcal{E})$, where
$[m]=\{1,2,\ldots,m\}$ is the set of nodes (agents) and $\mathcal{E}\subseteq [m]\times [m]$  is the edge set of ordered node pairs describing the interaction among agents.
For a nonnegative weighting matrix $W=\{w_{ij}\}\in\mathbb{R}^{m\times m}$, we define the induced directed graph as $\mathcal{G}_W=([m],\mathcal{E}_W)$, where
the directed edge $(i,j)$ from agent $j$ to agent $i$ exists,
 i.e., $(i,j)\in \mathcal{E}_W$ if and only if $w_{ij}>0$.
For an agent $i\in[m]$,
its in-neighbor set
$\mathbb{N}^{\rm in}_i$ is defined as the collection of agents $j$ such that $w_{ij}>0$; similarly,
the out-neighbor set $\mathbb{N}^{\rm out}_i$ of agent $i$ is the collection of agents $j$ such that $w_{ji}>0$.

The optimization problem (\ref{eq:optimization_formulation1}) can be
reformulated as the following equivalent multi-agent optimization
problem:
\vspace{-0.1cm}
\begin{equation}\label{eq:optimization_formulation2}
\vspace{-0.1cm}
\min\limits_{x\in\mathbb{R}^{md}}f(x)\triangleq
\frac{1}{m}\sum_{i=1}^m f_i(x_i)\:\: {\rm s.t.}\:\:
x_1=x_2=\cdots=x_m
\end{equation}
where $x_i\in\mathbb{R}^d$ is  agent $i$'s decision variable  and
 the collection
of the agents' variables is $x=[x_1^T, x_2^T, \ldots, x_m^T]^T\in\mathbb{R}^{md}$.

We make the following assumption on objective functions.
\begin{Assumption 2}\label{assumption:objective_functions}
Problem~\eqref{eq:optimization_formulation1}  has an optimal solution $\theta^{\ast}$.
The objective function   $F(\cdot)$ is convex and each $f_i(\cdot)$ has
Lipschitz continuous gradients over $\mathbb{R}^d$, i.e., for some
$L>0$,
\[\|\nabla f_i(u)-\nabla f_i(v)\|\le L \|u-v\|,\quad\forall i\in [m]  \:\:{\rm and}\:\:\forall u,v\in\mathbb{R}^d\]
\end{Assumption 2}

Under Assumption \ref{assumption:objective_functions}, the optimization problem~\eqref{eq:optimization_formulation2}  has an optimal
solution $ x^{\ast}=
[(\theta^{\ast})^T,(\theta^{\ast})^T,\ldots,(\theta^{\ast})^T]^T\in\mathbb{R}^{md}$.

 In the analysis of our methods, we use the following results.
\begin{Lemma 1}[\cite{polyak87}, Lemma 11, page 50]\label{Lemma-polyak}
Let $\{v^k\}$, $\{u^k\}$, $\{\alpha^k\}$, and $\{\b^k\}$ be random nonnegative scalar sequences such that
$\sum_{k=0}^\infty \a^k<\infty$ and $\sum_{k=0}^\infty \b^k<\infty$ \as
and
\[
\begin{aligned}
\mathbb{E}\left[v^{k+1}|\mathcal{F}^k \right]\le(1+\a^k) v^k -u^k+\b^k,\quad \forall k\geq 0 \quad \as
\end{aligned}
\]
where 
$\mathcal{F}^k=\{v^\ell,u^\ell,\a^\ell,\b^\ell;\, 0\le \ell\le k\}$.
Then $\sum_{k=0}^{\infty}u^k<\infty$ and $\lim_{k\to\infty}v^k=v$ for a random variable $v\geq 0$ \as
\end{Lemma 1}

\begin{Lemma 2}\label{Lemma-polyak_2}
Let $\{v^k\}$,$\{\a^k\}$, and $\{p^k\}$ be random nonnegative scalar sequences, and
$\{q^k\}$ be a deterministic nonnegative scalar sequence satisfying
$\sum_{k=0}^\infty \a^k<\infty$ {\it a.s.},
$\sum_{k=0}^\infty q^k=\infty$, $\sum_{k=0}^\infty p^k<\infty$ {\it a.s.},
and the following inequality
\[
\mathbb{E}\left[v^{k+1}|\mathcal{F}^k\right]\le(1+\a^k-q^k) v^k +p^k,\quad \forall k\geq 0\quad\as
\]
where $\mathcal{F}^k=\{v^\ell,\a^\ell,p^\ell; 0\le \ell\le k\}$.
Then, $\sum_{k=0}^{\infty}q^k v^k<\infty$ and
$\lim_{k\to\infty} v^k=0$ hold \as
\end{Lemma 2}
\begin{proof}
From the given relation we have \as
\be\label{eq-supmart}
\mathbb{E}\left[v^{k+1}|\mathcal{F}^k\right]\le (1+\a^k)v^k -q^k v^k +p^k,\quad \forall k\geq 0
\ee
By Lemma~\ref{Lemma-polyak} with $u^k=q^k v^k$, and $\b^k=p^k$, it follows that
$\sum_{k=0}^{\infty}q^k v^k<\infty$ and $\lim_{k\to\infty}v^k=v$ for a random variable $v\geq 0$ \as
Since $\sum_{k=0}^\infty q^k=\infty$, it follows that
$\liminf_{k\to\infty} v^k=0$ \as This and the fact  $v^k\to v$ \as imply that
$\lim_{k\to\infty}v^k=0$ \as
\end{proof}
%


\begin{Lemma 3}\label{lem-opt}
Consider the problem $\min_{z \in \R^d} \phi(z)$,
where $\phi:\mathbb{R}^d\to\mathbb{R}$ is a
continuous function. Assume that
the optimal solution set $Z^*$ of the problem is nonempty.
Let $\{z^k\}$ be a random sequence such that for any optimal solution $z^*\in Z^*$,
\[
\begin{aligned}
&\mathbb{E}\left[\|z^{k+1}-z^*\|^2|\mathcal{F}^k\right]\\
&\le (1+\a^k)\|z^k - z^*\|^2  - \eta^k\left (\phi(z^k) - \phi(z^*)\right) +\b^k,\: \forall k\geq 0
\end{aligned}
\]
holds {\it a.s.},
 where
$\mathcal{F}^k=\{z^\ell,\a^\ell,\b^\ell,\ \ell=0,1,\ldots,k\}$,
$\{\a^k\}$ and $\{\b^k\}$ are random nonnegative scalar sequences
satisfying $\sum_{k=0}^\infty \a^k<\infty$, $\sum_{k=0}^\infty \b^k<\infty$ \as,
while $\{\eta^k\}$ is a deterministic nonnegative scalar sequence with
$\sum_{k=0}^\infty \eta^k=\infty$.
Then, $\{z^k\}$ converges \as to some solution $z^*\in Z^*$.
\end{Lemma 3}
\begin{proof}
By letting $z=z^{\ast}$ for an arbitrary $z^{\ast}\in Z^{\ast}$ and
 defining  $\phi^{\ast}=\min_{z\in\mathbb{R}^m}\phi(z)$, we obtain \as for all $k$,
\[
\mathbb{E}\left[\|z^{k+1}\hspace{-0.1cm}-z^{\ast}\|^2|\mathcal{F}^k\right]\hspace{-0.07cm}\leq \hspace{-0.07cm} (1+\alpha^k)\|z_k-z^{\ast}\|
- \eta^k(\phi(z^k)-\phi^{\ast}) +\beta^k
\]
Thus, all the conditions of Lemma \ref{Lemma-polyak} are satisfied, yielding
\begin{equation}\label{eq:conv}
    \left\{\|z^k-z^{\ast}\|\right\} {\rm converges\: for \: each\:} z^{\ast}\in Z^{\ast}\quad a.s.
\end{equation}
\vspace{-0.4cm}
\begin{equation}\label{eq:summable}
\sum_{k=0}^{\infty} \eta^k(\phi(z^k)-\phi^{\ast})  <\infty \quad a.s.
\end{equation}
 From (\ref{eq:summable}) and $\sum_{k=0}^{\infty}\eta^k=\infty$ we have
$\liminf_{k\rightarrow \infty} \phi(z^k)=\phi^{\ast}$  \as
Let $\{z^{k_{\ell}}\}$ be a subsequence such that almost surely
\vspace{-0.1cm}
\begin{equation}\label{eq:converge_subsequence}
 \lim_{\ell\rightarrow\infty}\phi(z^{k_{\ell}})=\liminf_{k\rightarrow \infty}\phi(z^k)=\phi^{\ast}
\vspace{-0.1cm}
\end{equation}
 Relation (\ref{eq:conv}) implies that the sequence $\{z^k\}$ is bounded \as Thus,
we can assume without loss of generality that $\{z^{k_{\ell}}\}$ converges {\it a.s.} to some
$\tilde{z}$ (for otherwise, we can in turn select a convergent subsequence of $\{z^{k_{\ell}}\}$). Therefore, by the continuity of  $\phi$, one has $\lim_{\ell\rightarrow\infty}\phi(z^{k_{\ell}})=\phi(\tilde{z})$  {\it a.s.}, which in combination with (\ref{eq:converge_subsequence}) implies that $\tilde{z}\in Z^{\ast}$ \as By letting $z^{\ast}=\tilde{z}$ in (\ref{eq:conv}), we see that $z^k$ converges to $\tilde{z}$ \as
\end{proof}

\begin{Lemma 4}\label{le:chung}
Let $\{v^k\}$ be a nonnegative sequence, and  $\{\a^k\}$  and $\{\b^k\}$ be positive sequences satisfying $\sum_{k=0}^{\infty}\a^k=\infty$, $\lim_{k\rightarrow \infty} \a^k =0$, and $\lim_{k\rightarrow \infty}\frac{\beta^k}{\alpha^k}\rightarrow0$ with a polynomial  decay rate. If there exists a $K\geq 0$ such that $ v^{k+1} \le(1-\a^k) v^k +\b^k$ holds for all $k\geq K$,
 then we always have $v^k\leq C \frac{\b^k}{\a^k}$  for all $k$, where $C$ is some constant.
\end{Lemma 4}
\begin{proof}
  The derivation follows the same line of reasoning in Lemma 4 of \cite{chung1954stochastic} and is omitted here.
\end{proof}
{
\subsection{On differential privacy}
We consider Laplace   noise for DP. For a constant $\nu>0$, ${\rm Lap}(\nu)$ denotes the Laplace distribution with probability density function $\frac{1}{2\nu}e^{-\frac{|x|}{\nu}}$. This distribution has mean zero and variance $2\nu^2$.
Following \cite{Huang15}, for the convenience of DP analysis, we represent the distributed optimization problem $\mathcal{P}$ in (\ref{eq:optimization_formulation1}) by four parameters ($\mathcal{X},\mathcal{S},F,\mathcal{G}_W$), where
 $\mathcal{X}=\mathbb{R}^n$ is the domain of optimization,  $\mathcal{S}\subseteq\{\mathbb{R}^n\mapsto\mathbb{R}\}$ is a set of real-valued objective functions, with $f_i\in\mathcal{S}$, and $F(x)\triangleq\frac{1}{m}\sum_{i=1}^m f_i(x)$, and $\mathcal{G}_W$ is the induced graph by matrix $W$. Then we define adjacency as follows:

\begin{Definition 1}\label{de:adjacency}
Two distributed optimization problems $\mathcal{P}$ and $\mathcal{P}'$ are adjacent if the following conditions hold:
\begin{itemize}
\item $\mathcal{X}=\mathcal{X}'$, $\mathcal{S}=\mathcal{S}'$, and $\mathcal{G}_W=\mathcal{G}_W'$, i.e., the domain of optimization, the set of individual objective functions, and the communication graphs are identical;
\item there exists an $i\in[m]$ such that $f_i\neq f_i'$ but $f_j=f_j'$ for all $j\in[m],\,j\neq i$;
\item  the different objective functions $f_i$ and $f'_i$  have similar behaviors   around $\theta^\ast$, the  solution of $\mathcal{P}$. More specifically,  there exits some $\delta>0$  such that for all $v$ and $v'$ in  $B_\delta(\theta^\ast)\triangleq\{u:u\in\mathbb{R}^d, \|u-\theta^\ast\|<\delta\}$, we have $\nabla f_i(v)=\nabla f'_i(v')$.
\end{itemize}
\end{Definition 1}

It can be seen that two distributed optimization problems are adjacent if and only if one agent changes its individual objective function while all others parameters are identical.

\begin{Remark 1}
 In Definition \ref{de:adjacency}, since  the change of an objective function from $f_i$ to $f'_i$ in the second condition  can be arbitrary,
 additional restrictions have to imposed to ensure rigorous DP in distributed optimization.  Different from \cite{Huang15} which restricts all gradients  to be uniformly bounded,   we add the third condition, which, as shown later, allows us to ensure  rigorous DP while maintaining provable convergence to the optimal solution.
\end{Remark 1}

 Given a distributed optimization problem, we represent an execution of such an algorithm as $\mathcal{A}$, which is an infinite sequence of the optimization variables, i.e., $\mathcal{A}=\{x^0,x^1,\cdots\}$. We consider adversaries that can observe all communicated messages in the network. Therefore, the observation part of an execution is the infinite sequence of shared messages, which is represented by $\mathcal{O}$. Given a distributed optimization problem $\mathcal{P}$ and an initial state $x^0$, we define the observation mapping as {$\mathcal{R}_{\mathcal{P},x^0}(\mathcal{A})\triangleq \mathcal{O}$}. Given a distributed optimization problem $\mathcal{P}$, observation sequence $\mathcal{O}$, and an initial state $x^0$,  {$\mathcal{R}_{\mathcal{P},x^0}^{-1}(\mathcal{O})$} is the set of executions $\mathcal{A}$ that can generate  observation $\mathcal{O}$.
 \begin{Definition 1}
   ($\epsilon$-DP \cite{Huang15}). For a given $\epsilon>0$, an iterative   algorithm for problem~(\ref{eq:optimization_formulation1}) is $\epsilon$-differentially private if for any two adjacent $\mathcal{P}$ and $\mathcal{P}'$, any set of observation sequences $\mathcal{O}_s\subseteq\mathbb{O}$ (with $\mathbb{O}$ denoting the set of all possible observation sequences), and any initial state ${x}^0$, we always have
    \vspace{-2mm}
    \begin{equation}
        {\mathbb{P}[\mathcal{R}_{\mathcal{P},x^0}^{-1}\left(\mathcal{O}_s\right)]}\leq e^\epsilon{\mathbb{P}[\mathcal{R}_{\mathcal{P}',x^0}^{-1}\left(\mathcal{O}_s\right)]}
    \end{equation}
    where the probability $\mathbb{P}$ is taken over the randomness over iteration processes.
 \end{Definition 1}

The definition of $\epsilon$-DP ensures that an adversary having access to all shared messages in the network cannot gain information with a  significant probability of any participating agent's objective function. It can also be seen that a smaller $\epsilon$ means a higher level of privacy protection.
}
\section{
Static-consensus gradient methods  for differentially-private distributed optimization}\label{sec-static}
In this section, we  tailor a static-consensus based distributed gradient method
to construct  a differentially-private distributed method with almost sure convergence to
 an  optimal solution.
 The agent interaction strength is captured by a weight matrix
$W=\{w_{ij}\}$, where $w_{ij}>0$ if there is a link from agent  $j$ to agent $i$,
and $w_{ij}=0$ otherwise.
We let $w_{ii}\triangleq-\sum_{j\in\mathbb{N}^{\rm in}_i}w_{ij}$  for all $i\in [m]$,
where $\mathbb{N}^{\rm in}_i$ is the in-neighbor set of agent $i$.
We make the following assumption on $W$:

\begin{Assumption 1}\label{assumption:W}
 The matrix  $W=\{w_{ij}\}\in \mathbb{R}^{m\times m}$ is symmetric and satisfies
    ${\bf 1}^TW={\bf
  0}^T$, $W{\bf 1}={\bf
  0}$, $ \|I+W-\frac{{\bf 1}{\bf 1}^T}{m}\|<1$.
\end{Assumption 1}

Assumption~\ref{assumption:W} ensures that the interaction graph induced by $W$ is balanced and  connected, i.e., there is a  path
from each agent to every other agent.

To achieve a strong DP, independent DP-noise should be injected  repeatedly in every round of message sharing and, hence, constantly affects the algorithm
through inter-agent interactions, leading to significant reduction in optimization accuracy.
Motivated by this observation, we propose to gradually weaken inter-agent interactions to reduce the influence of DP-noise on optimization accuracy. Interestingly, we prove that by judiciously designing the interaction weakening mechanism, we can still ensure  convergence of all agents   to a common optimal solution  even in the presence of persistent
DP-noise.


\noindent\rule{0.49\textwidth}{0.5pt}
\noindent\textbf{Algorithm 1: DP-oriented static-consensus based distributed optimization}

\noindent\rule{0.49\textwidth}{0.5pt}
\begin{enumerate}[wide, labelwidth=!, labelindent=0pt]
    \item[] Parameters: Stepsize $\lambda^k$ and
    weakening factor $\gamma^k$.
    \item[] Every agent $i$ maintains one state     $x_i^k$, which is initialized with a random vector in $\mathbb{R}^d$.
    \item[] {\bf for  $k=1,2,\ldots$ do}
    \begin{enumerate}
        \item Every agent $j$ adds persistent DP-noise   $\zeta_j^{k}$ 
        to its state
    $x_j^k$,  and then sends the obscured state $x_j^k+\zeta_j^{k}$ to agent
        $i\in\mathbb{N}_j^{\rm out}$.
        \item After receiving  $x_j^k+\zeta_j^k$ from all $j\in\mathbb{N}_i^{\rm in}$, agent $i$ updates its state  as follows:
        \begin{equation}\label{eq:diminishing_update}
        \hspace{-0.5cm}
        \begin{aligned}
             x_i^{k+1}\hspace{-0.1cm}&=\hspace{-0.1cm}x_i^k+\sum_{j\in \mathbb{N}_i^{\rm
            in}}\gamma^k w_{ij}(x_j^k+\zeta_j^k-x_i^k)-\lambda^k\nabla f_i(x_i^k)
        \end{aligned}
        \end{equation}
                \item {\bf end}
    \end{enumerate}
\end{enumerate}
\vspace{-0.1cm} \rule{0.49\textwidth}{0.5pt}

 The sequence $\{\gamma^k\}$ diminishes with time and is used to suppress the influence of persistent DP-noise $\zeta_j^k$ on the convergence point of the iterates.
The stepsize sequence $\{\lambda^k\}$ and attenuation sequence $\{\gamma^k\}$
have to be designed appropriately to guarantee the almost sure convergence of all $\{x_i^k\}$   to a common optimal solution $\theta^{\ast}$.
The persistent DP-noise processes $\{\zeta_i^k\}, i\in[m]$ have zero-mean and
$\gamma^k$-bounded  (conditional) variances, to be specified later in Assumption \ref{assumption:dp-noise}.

\subsection{Convergence analysis}
We   have to extend Lemma~\ref{Lemma-polyak} to deal with random vectors.

\begin{Lemma 1}\label{th-dsystem}
Let  $\{\bv^k\}\subset \mathbb{R}^d$
and $\{\bu^k\}\subset \mathbb{R}^p$ be random nonnegative
vector sequences, and $\{a^k\}$ and $\{b^k\}$ be random nonnegative scalar sequences   such that
\[
\mathbb{E}\left[\bv^{k+1}|\mathcal{F}^k\right]\le (V^k+a^k{\bf 1}{\bf1}^T)\bv^k +b^k{\bf 1} -H^k\bu^k,\quad \forall k\geq 0
\]
holds {\it a.s.}, where $\{V^k\}$ and $\{H^k\}$ are random sequences of
nonnegative matrices and
$\mathbb{E}\left[\bv^{k+1}|\mathcal{F}^k \right]$ denotes the conditional expectation given
 $\bv^\ell,\bu^\ell,a^\ell,b^\ell,V^\ell,H^\ell$ for $\ell=0,1,\ldots,k$.
Assume that $\{a^k\}$ and $\{b^k\}$ satisfy
$\sum_{k=0}^\infty a^k<\infty$ and $\sum_{k=0}^\infty b^k<\infty$ {\it a.s.}, and
that there exists a (deterministic) vector $\pi>0$ such that
$\pi^T V^k\le \pi^T$ and $\pi^TH^k\ge 0$ hold {\it a.s.} for all $k\geq 0$.
Then, we have
1) $\{\pi^T\bv^k\}$ converges to some random variable $\pi^T\bv\geq 0$ {\it a.s.}; 2) $\{\bv^k\}$ is bounded {\it a.s.}, and
3) $\sum_{ k=0 }^\infty \pi^TH^k\bu^k<\infty$ holds {\it a.s.}
\end{Lemma 1}

\begin{proof}
By multiplying the given relation for $\bv^{k+1}$ with $\pi$ and using $\pi^TV^k\le\pi^T$ and the nonnegativity of  $\bv^k$,
we obtain
\[
\mathbb{E}\hspace{-0.07cm}\left[\pi^T\bv^{k+1}|\mathcal{F}^k\right]\hspace{-0.07cm}\le\hspace{-0.07cm} \pi^T\bv^k\hspace{-0.07cm} +a^k(\pi^T {\bf 1})({\bf1}^T\bv^k) +b^k\pi^T{\bf 1}-\pi^TH^k\bu^k
\]
Since $\pi>0$, we have $\pi_{\min}=\min_i\{\pi_i\}>0$, which yields
${\bf1}^T\bv^k=\frac{1}{\pi_{\min}} \,\pi_{\min}{\bf1}^T\bv^k\le
\frac{1}{\pi_{\min}} \,\pi^T\bv^k$,
where the inequality holds since $\bv^k\ge0$.
So, one obtains
\[\mathbb{E}\left[\pi^T\bv^{k+1}|\mathcal{F}^k\right]\le \left(1+a^k\frac{\pi^T {\bf 1}}{\pi_{\min}}\right)\pi^T\bv^k +b^k\pi^T{\bf 1} -\pi^TH^k\bu^k\]
By our assumption,  $\pi^TH^k\bu^k\ge0$ holds for all $k$ {\it a.s.}
Thus, the preceding relation implies that
the conditions of Lemma~\ref{Lemma-polyak} are satisfied
with $v^k = \pi^T\bv^k$, $\a^k = a^k\pi^T {\bf 1}/\pi_{\min}$ and $\b^k = b^k \pi^T{\bf 1}$.
So by Lemma~\ref{Lemma-polyak},
it follows that  $\lim_{k\to\infty}\pi^T\bv^k$ exists {\it a.s.}
Consequently, $\{\pi^T\bv^k\}$ is bounded {\it a.s.}, and
since $\{\bv^k\}$ is nonnegative  and $\pi>0$, it follows that $\{\bv^k\}$ is also bounded {\it a.s.}
By Lemma \ref{Lemma-polyak}, we have
$\sum_{k=0}^\infty \pi^TH^k\bu^k<\infty$ {\it a.s.}
\end{proof}

Based on Lemma \ref{lem-opt} and Lemma \ref{th-dsystem}, we can prove the following
general convergence results for static-consensus based distributed algorithms for problem (\ref{eq:optimization_formulation1}).

\begin{Lemma 2}\label{th-main_decreasing}
Assume that problem (1) has a solution. Suppose that a distributed algorithm generates sequences
$\{x_i^k\}\subseteq\R^d$ such that
\as we have for any optimal solution $\theta^*$,
\begin{equation}\label{eq:Theorem_decreasing}
\begin{aligned}
&\left[\begin{array}{c}
\mathbb{E}\left[\|\bar x^{k+1}-\theta^*\|^2|\mathcal{F}^k\right]\cr
\mathbb{E}\left[\sum_{i=1}^m\|x_i^{k+1}-\bar x^{k+1}\|^2|\mathcal{F}^k\right]\end{array}
\right]
\\
&\le \left( \left[\begin{array}{cc}
1 & \frac{\gamma^k}{m}\cr
0& 1-\kappa\gamma^k\cr
\end{array}\right]
+a^k {\bf 1}{\bf 1}^T\right)\left[\begin{array}{c}\|\bar x^k-\theta^*\|^2\cr
\sum_{i=1}^m\|x_i^k-\bar x^k\|^2\end{array}\right]&&\cr
&\quad+b^k{\bf 1} - c^k \left[\begin{array}{c}
 F(\bar x^k)- F(\theta^*) \cr
 0\end{array}\right],\quad\forall k\geq 0
 \end{aligned}
\end{equation}
where  $\bar{x}^k=\frac{1}{m}\sum_{i=1}^m x_i^k$,
$\mathcal{F}^k=\{x_i^\ell, \, i\in[m],\, 0\le \ell\le k\}$,
the random
nonnegative scalar sequences $\{a^k\}$, $\{b^k\}$ satisfy $\sum_{k=0}^\infty a^k<\infty$ and $\sum_{k=0}^\infty b^k<\infty$ {\it a.s.},  the deterministic  nonnegative sequences $\{c^k\}$ and $\{\gamma^k\}$  satisfy
$
\sum_{k=0}^\infty c^k=\infty$ and $
\sum_{k=0}^\infty \gamma^k=\infty
$, and  the scalar $\kappa>0$ satisfies  $\kappa\gamma^k<1$ for all $k\geq 0$.
Then, we have
$\lim_{k\to\infty}\|x_i^k - \bar x^k\|=0$ \as for all   $i$,
and there is a solution $\tilde\theta^*$ such that
$\lim_{k\to\infty}\|\bar x^k-\tilde\theta^*\|=0$ \as
\end{Lemma 2}
\begin{proof}
  See Appendix A.
\end{proof}

Using Lemma~ \ref{th-main_decreasing}, we are in position to
 establish convergence  of Algorithm 1 assuming that persistent DP-noise satisfies the following assumption.
\begin{Assumption 1}\label{assumption:dp-noise}
For every $i\in[m]$ and every $k$, conditional on the state $x_i^k$,
the random noise $\zeta_i^k$ satisfies $
\mathbb{E}\left[\zeta_i^k\mid x_i^k\right]=0$ and $\mathbb{E}\left[\|\zeta_i^k\|^2\mid x_i^k\right]=(\sigma_{i}^k)^2$ for all  $k\ge0$, and
\begin{equation}\label{eq:condition_assumption1}
\sum_{k=0}^\infty (\gamma^k)^2\, \max_{i\in[m]}(\sigma_{i}^k)^2 <\infty
\end{equation} where $\{\gamma^k\}$ is the attenuation sequence from Algorithm 1.
The initial random vectors satisfy
$\mathbb{E}\left[\|x_i^0\|^2\right]<\infty$,  $\forall i\in[m]$.
\end{Assumption 1}
{
\begin{Remark 1}
Given that $\gamma^k$ decreases with time, (\ref{eq:condition_assumption1}) can be satisfied even when $\{\sigma_i^k\}$ increases with time. For example, under $\gamma^k=\mathcal{O}(\frac{1}{k^{0.9}})$, an increasing $\{\sigma_i^k\}$ with increasing rate no faster than $\mathcal{O}(k^{0.3})$
still satisfies the summable condition in (\ref{eq:condition_assumption1}). Allowing $\{\sigma_i^k\}$ to  increase with time is key to enabling the strong $\epsilon$-DP in Theorem \ref{th:DP_Algorithm1}.
\end{Remark 1}}

\begin{Theorem 1}\label{theorem:convergence_algorithm_1}
Under Assumption \ref{assumption:objective_functions}, Assumption \ref{assumption:W},
 and Assumption~\ref{assumption:dp-noise},
 Algorithm~1 converges to a solution of problem~(1)  {\it a.s.}
 when nonnegative sequences $\{\gamma^k\}$ and $\{\lambda^k\}$ satisfy $
\sum_{k=0}^\infty \gamma^k=\infty$, $\sum_{k=0}^\infty \lambda^k=\infty$,  and $\sum_{k=0}^\infty \frac{(\lambda^k)^2}{\gamma^k}<\infty$.
\end{Theorem 1}

\begin{proof}
See Appendix B.
\end{proof}

\begin{Remark 1}
Communication imperfections  can be modeled as channel noises \cite{srivastava2011distributed,kar2008distributed}, which can be regarded as the DP-noise   here.  Therefore, Algorithm 1 can also  counteract such communication imperfections in distributed optimization.
\end{Remark 1}

{

\begin{Remark 1}\label{re:convergence_algorithm1}
Because the evolution of $x_i^k$ to the optimal solution satisfies the conditions in Lemma \ref{th-main_decreasing}, we can leverage Lemma \ref{th-main_decreasing} to examine the convergence speed.
 From Lemma \ref{le:chung}, the relationship in (\ref{eq:speed_in_lemma5})   implies that   $\sum_{i=1}^m \|x_i^k - \bar x^k\|^2$  decreases to zero   no slower than $\mathcal{O}((\frac{\lambda^k}{\gamma^k})^2)$, and hence we have $x_i^k$ converging to $\bar{x}^k$  no slower than $\mathcal{O}(\frac{\lambda^k}{\gamma^k})$ (note $\beta^k$ is on the order of $\frac{(\lambda^k)^2}{\gamma^k}$ from the proof of Theorem \ref{theorem:convergence_algorithm_1}). Moreover, when $F$ is strongly convex,  (\ref{eq:x_opt}) implies that $\bar{x}^k$ converges to $\theta^\ast$   no slower than $\mathcal{O}((\frac{\lambda^k}{\gamma^k})^{0.5})$ using Lemma \ref{le:chung}. Therefore, the convergence of every $x_i^k$ to $\theta^\ast$, which is equivalent to the combination of the convergence of $x_i^k$ to $\bar{x}^k$ and the convergence of $\bar{x}^k$ to $\theta^\ast$,  should be no slower than $\mathcal{O}((\frac{\lambda^k}{\gamma^k})^{0.5})$.    Moreover, from the proof of the theorem, it can be seen that the decreasing speed of  $\|x^{k}-{\bf 1}\otimes \bar{x}^{k}\|^2 $ (where $\otimes$ is the Kronecker product) increases with an increase in $|\nu|$, which corresponds to the spectral radius of $W$. Therefore, the decreasing speed of  $\|x^{k}-{\bf 1}\otimes \bar{x}^{k}\|^2 $ to zero increases with an increase in the spectral radius of $W$ defined in Assumption \ref{assumption:W}.
\end{Remark 1}

\subsection{Privacy analysis}\label{se:privacy_Algorithm1}
Similar to \cite{Huang15}, 
we  define the sensitivity of an algorithm to problem (\ref{eq:optimization_formulation1}) as follows:
\begin{Definition 1}\label{de:sensitivity}
  At each iteration $k$, any initial state $x^0$ and any adjacent distributed optimization problems $\mathcal{P}$ and $\mathcal{P'}$,  the sensitivity of an algorithm is
  \begin{equation}
  \Delta^k\triangleq \sup\limits_{\mathcal{O}\in\mathbb{O}}\left\{\sup\limits_{x\in\mathcal{R}_{\mathcal{P},x^0}^{-1}(\mathcal{O}),\:x'\in\mathcal{R}_{\mathcal{P'},x^0}^{-1}(\mathcal{O})}\|x^k-x'^k\|_1\right\}
  \end{equation}
\end{Definition 1}

\begin{Lemma 1}\label{Le:Laplacian}
At each iteration $k$, if each agent adds a noise vector $\zeta_i^k\in\mathbb{R}^p$ consisting of $p$ independent Laplace noises with  parameter $\nu^k$ such that $\sum_{k=1}^{T}\frac{\Delta^k}{\nu^k}\leq \epsilon$, then  Algorithm 1 is $\epsilon$-differentially private for iterations from $k=0$ to $k=T$.
\end{Lemma 1}
\begin{proof}
The lemma can be obtained following the same line of reasoning of Lemma 2 in  \cite{Huang15}.
\end{proof}


\begin{Theorem 1}\label{th:DP_Algorithm1}
Under Assumptions \ref{assumption:objective_functions} and \ref{assumption:W}, if nonnegative sequences $\{\lambda^k\}$ and $\{\gamma^k\}$ satisfy the conditions in Theorem \ref{theorem:convergence_algorithm_1}, and all elements of $\zeta_i^k$ are drawn independently from  Laplace distribution ${\rm Lap}(\nu^k)$ with $(\sigma_i^k)^2=2(\nu^k)^2$ satisfying Assumption \ref{assumption:dp-noise}, then all agents in Algorithm 1 will converge \as to an optimal solution. Moreover, 
\begin{enumerate}
\item For any finite number of iterations $T$, Algorithm 1 is  $\epsilon$-differentially private with the cumulative privacy budget bounded by $\epsilon\leq \sum_{k=1}^{T}\frac{C\varsigma^k}{\nu^k}$   where $\varsigma^k\triangleq \sum_{p=1}^{k-1}(\Pi_{q=p}^{k-1}(1-\bar{w}\gamma^{q}))\lambda^{p-1}+\lambda^{k-1}$, $\bar{w}\triangleq\min_i\{|w_{ii}|\}$, and $C\triangleq \max_{i\in[m],0\leq k\leq T-1}\{\|\nabla f_i(x_i^k)-\nabla f'_i({x'_i}^k)\|_1\}$ (note that $C$ is always finite since the algorithm ensures convergence in both $\mathcal{P}$ and $\mathcal{P}'$);
\item  The cumulative privacy budget is  finite for $T\rightarrow\infty$  when the sequence  $\{\frac{\lambda^k}{\nu^k}\}$ is summable.
\end{enumerate}
\end{Theorem 1}
\begin{proof}
Since the Laplace noise satisfies Assumption \ref{assumption:dp-noise}, the convergence results follow naturally from    Theorem \ref{theorem:convergence_algorithm_1}.

To prove the  statements on privacy, we first analyze the sensitivity of Algorithm 1.
 Given two adjacent distributed optimization problems $\mathcal{P}$ and $\mathcal{P'}$, for any given fixed observation $\mathcal{O}$ and initial state $x^0$, the sensitivity depends on $\|x^{k}-x'^{k}\|_1$ according to Definition \ref{de:sensitivity}. Since in $\mathcal{P}$ and $\mathcal{P'}$, there is only one objective function that is different, we  represent this different objective function as the  $i$th one, i.e., $f_i$ in $\mathcal{P}$ and $f'_i$ in $\mathcal{P}'$, without loss of generality.

 Because the initial conditions, objective functions, and observations of $\mathcal{P}$ and $\mathcal{P'}$  are identical for $j\neq i$, we have $x_j^k={x'_j}^k$ for all $j\neq i$ and $k$. Therefore, $\|x^{k}-x'^{k}\|_1$ is always equal to $\|x_i^{k}-{x'_i}^{k}\|_1$.

 According to   Algorithm 1, we can arrive at
 \[
 \begin{aligned}
 &x_i^{k+1}-{x'_i}^{k+1}&=(1+w_{ii}\gamma^k)(x_i^k-{x'_i}^k)-\lambda^k(g_i^k-{g'_i}^k),
 \end{aligned}
 \]
 where we have represented $\nabla f_i(x_i^k)$ and $\nabla f'_i({x'_i}^k)$ as $g_i^k$ and ${g'_i}^k$, respectively,  for notational simplicity.
 Note that we have also used the definition $w_{ii}\triangleq-\sum_{j\in\mathbb{N}_i}w_{ij}$ and  the fact that the observations $x_j^k+\zeta_j^k$ and ${x'_j}^k+{\zeta'_j}^k$ are the same.

 Hence, the sensitivity $\Delta^k$ satisfies
 \[
 \Delta^{k+1}\leq (1-|w_{ii}|\gamma^k)\Delta^{k}+\lambda^k \|g_i^k-{g'_i}^k\|_1.
 \]
 which, implies the first statement by iteration using Lemma \ref{Le:Laplacian}.

 For the infinity horizon result in the second statement, we exploit the fact that   our algorithm ensures convergence in both $\mathcal{P}$ and $\mathcal{P'}$. This means that  $\|g_i^k-{g'_i}^k\|_1=0$ will  be satisfied when $k$ is large enough using the third condition in Definition \ref{de:adjacency} (see Remark \ref{re:convergence_algorithm1} for convergence rate analysis). Furthermore, the ensured convergence also means that $\|g_i^k-{g'_i}^k\|_1$ is always  bounded. Hence, there always exists some constant ${C}$ such that the sequence $\{\|g_i^k-{g'_i}^k\|_1\}$ is upper bounded by the sequence $\{ {C}\gamma^k\}$.

 Therefore, according to Lemma \ref{le:chung},   there always exists a constant $\bar{C}$ such that  $\Delta^k\leq \bar{C}\lambda^k$ holds.
Using Lemma \ref{Le:Laplacian}, we can easily obtain $\epsilon\leq \sum_{k=1}^{T}\frac{\bar{C}\lambda^k}{\nu^k}$. Hence, $\epsilon$ will  be finite even when $T$ tends to infinity if  the sequence $\{\frac{\lambda^k}{\nu^k}\}$ is summable, i.e.,  $\sum_{k=0}^{\infty}\frac{\lambda^k}{\nu^k}<\infty$.
\end{proof}

Different from \cite{Huang15}  which has to use a summable stepsize  (geometrically-decreasing stepsize, more specifically) to ensure a finite   privacy budget $\epsilon$ when $k\rightarrow\infty$, here we ensure   a finite  $\epsilon$  even when the stepsize sequence is non-summable. Allowing stepsize sequences to be non-summable is key to avoiding optimization errors in \cite{Huang15} and achieve almost sure convergence. In fact, to our knowledge, this is the first time that almost-sure convergence is achieved under rigorous $\epsilon$-DP for an infinite number of iterations.

\begin{Remark 1}
  In Theorem \ref{th:DP_Algorithm1}, to ensure that the privacy budget   is finite even when  $k\rightarrow \infty$, the Laplace noise parameter $\nu^k$ has to increase with time since $\{\lambda^k\}$ is non-summable.  An increasing $\nu^k$ will make the relative level between noise $\zeta_i^k$ and signal $x_i^k$ increase  with time. However, since the increase  in $\nu^k$ is outweighed by the decrease of $\gamma^k$ (see  Assumption \ref{assumption:dp-noise}), the actual noise fed into the algorithm, i.e., $\gamma^k{\rm Lap}(\nu^k)$, still decays with time, which makes it possible for  Algorithm 1 to ensure  \as convergence to an optimal solution. Moreover, according to Theorem \ref{theorem:convergence_algorithm_1}, such \as convergence is not affected by scaling  $\nu^k$ by any constant coefficient $\frac{1}{\epsilon}>0$ so as to achieve any desired level of $\epsilon$-DP, as long as the Laplace noise parameter $\nu^k$ (with associated variance $(\sigma_i^k)^2=2(\nu^k)^2$) satisfies Assumption \ref{assumption:dp-noise}. 
\end{Remark 1}

}

\section{
gradient-tracking based methods for differentially private distributed optimization}\label{sec-dynamic}
In this section,  we  propose a DP-oriented gradient-tracking based distributed algorithm for general directed graphs and prove that it can ensure convergence to an optimal solution even under persistent DP-noise.
 In  gradient-tracking based algorithms, every agent $i\in[m]$ maintains and updates two iterates, $x_i^k$ and $y_i^k$, where  $y_i^k$ is an estimate of the ``joint agent" descent direction.
These two iterates are exchanged with local neighbors in two different communication networks, namely,  $\mathcal{G}_R$ and $\mathcal{G}_C$, which are, respectively, induced by matrices $R\in\R^{m\times m}$ and $C\in\R^{m\times m}$;
that is $(i,j)$ is a directed link in the graph $\mathcal{G}_R$ if and only if $R_{ij}>0$ and, similarly, $(i,j)$ is a directed link in $\mathcal{G}_C$ if and only if $C_{ij}>0$.
We make the following assumption on  $R$ and $C$. Note that,
$\mathcal{R}_{A^T}$ is identical to
$\mathcal{R}_{A}$ with the directions of edges reversed.

\begin{Assumption 4}\label{Assumption:push_pull topology}
The matrices $R,C\in\R^{m\times m}$ have nonnegative off-diagonal entries
($R_{ij}\geq 0$ and $C_{ij}\geq 0$ for all $i\neq j$). The induced graphs $\mathcal{G}_R$ and
$\mathcal{G}_{C^T}$ satisfy
\begin{enumerate}
  \item  $\mathcal{G}_R$ and $\mathcal{R}_{C^T}$ each contain at least one spanning tree;
  \item  There exists at least one node that is a root of spanning trees for both $\mathcal{G}_R$ and $\mathcal{R}_{C^T}$.
\end{enumerate}
\end{Assumption 4}

\begin{Remark 1}
The assumption on $\mathcal{G}_R$ and $\mathcal{R}_{C^T}$
 is weaker than requiring that both induced graphs of $R$ and $C$ to be strongly connected, which is assumed in most of the existing works.
\end{Remark 1}

\noindent\rule{0.49\textwidth}{0.5pt}
\noindent\textbf{Algorithm 2: DP-oriented gradient-tracking
based distributed optimization}

\vspace{-0.2cm}\noindent\rule{0.49\textwidth}{0.5pt}
\begin{enumerate}[wide, labelwidth=!, labelindent=0pt]
    \item[] Parameters: Stepsizes $\lambda^k$, $\alpha^k$ and weakening factors $\gamma_1^k$, $\gamma_2^k$.
    \item[] Every agent $i$ maintains two states  $x_i^k$ and
    $y_i^k$, which are initialized with a random point  $x_i^0\in\mathbb{R}^d$ and $y_i^0=\nabla f_i(x_i^0)$.
    \item[] {\bf for  $k=1,2,\cdots$ do}
    \begin{enumerate}
        \item Every agent $j$ injects zero-mean DP-noises $\zeta_j^k$ and $\xi_j^k$  to its states      $x_j^k$ and $y_j^k$, respectively.
        \item Agent $i$ pushes  $C_{li}(y_i^k+\xi_i^k)$ to each agent
        $l\in\mathbb{N}_{C,i}^{\rm out}$, and it pulls $x_j^k+\zeta_j^k$ from each $j\in\mathbb{N}_{R,i}^{\rm in}$, where the subscript  $R$ or $C$ in neighbor sets indicates the neighbors with respect to  the graphs induced by these matrices.
         \item Agent $i$ chooses $\gamma_1^k>0$ and $\gamma_2^k>0$ satisfying
        $1+\gamma_1^kR_{ii}>0$ and $1+\gamma_2^kC_{ii}>0$ with        \begin{equation}\label{eq:diagonal_entries}
                  R_{ii}=-\sum_{j\in\mathbb{N}_{R,i}^{\rm in}}R_{ij} ,\quad C_{ii}=-\sum_{j\in\mathbb{N}_{C,i}^{\rm out}}C_{ji}
                \end{equation}
        Then, agent $i$ updates its states as follows:
        \begin{equation}\label{eq:update}
        \begin{aligned}
            x_i^{k+1}&=(1+\gamma_1^k R_{ii})x_i^k+\gamma_1^k\sum_{j\in \mathbb{N}_{R,i}^{\rm
            in}}R_{ij}(x_j^k+\zeta_j^k)-\lambda^ky_i^{k}\\
            y_i^{k+1}&=(1-\alpha^k+\gamma^k_2C_{ii})y_i^k+\gamma_2^k\sum_{j\in \mathbb{N}_{C,i}^{\rm
            in}}C_{ij}(y_j^k+\xi_j^k)\\
            &\qquad+\nabla f_i(x_i^{k+1})-(1-\alpha^k)\nabla f_i(x_i^k)
        \end{aligned}
         \end{equation}

                  \item {\bf end}
    \end{enumerate}
\end{enumerate}
\vspace{-0.3cm} \rule{0.49\textwidth}{0.5pt}

 Note that  the definition of $R_{ii}$ and $C_{ii}$ in (\ref{eq:diagonal_entries}) ensures that   $R=\{ R_{ij}\}$ has zero row sums  and
$C=\{C_{ij}\}$ has zero column sums.

\subsection{Convergence analysis}
We will prove that, when the two sequences $\{\gamma_1^k\}$ and $\{\gamma_2^k\}$ are designed appropriately, all agents' $x$-iterates generated by Algorithm~2 converge to an optimal solution {\it a.s.}, as long as the  injected noises $\zeta_j^k$ and $\xi_j^k$
have zero-mean and $\gamma_1^k(\gamma_2^k)$
bounded    variances, to be specified later in Assumption \ref{assumption:dp-noises-intrack}.
 To this end, we first  extend Lemma \ref{Lemma-polyak_2} to vectors.

\begin{Lemma 1}\label{Theo:convergence_to_zero}
Let $\{\bv^k\}\subset \mathbb{R}^d$  be a sequence of non-negative random vectors and
 $\{b^k\}$ be a sequence of nonnegative random scalars such that
$\sum_{k=0}^\infty b^k<\infty$ {\it a.s.}
and
\[
\mathbb{E}\left[
\bv^{k+1}|\mathcal{F}^k\right]\le V^k\bv^k +b^k{\bf 1}, \quad \forall k\ge0\quad \as
\]
where $\{V^k\}$ is a sequence of
non-negative matrices  and $\mathcal{F}^k=\{\bv^\ell,b^\ell;0\le \ell\le k\}$.
Assume that there exist a vector $\pi>0$ and a deterministic scalar sequence $\{a^k\}$ satisfying
$a^k\in(0,1)$, $\sum_{k=0}^\infty a^k=\infty$,
and
$\pi^T V^k\le (1-a^k) \pi^T$  for all $k\ge0$.
Then, we have $\lim_{k\to\infty}\bv^k=0$ \as
\end{Lemma 1}
\begin{proof}
By multiplying the given relation for $\bv^{k+1}$ with $\pi$ and using
$\pi^TV^k\le(1-a^k)\pi^T$, we obtain the following relation due to
the nonnegativity of the vectors $\bv^k$:
\[
\mathbb{E}\left[\pi^T\bv^{k+1}|\mathcal{F}^k\right]\le (1-a^k)\pi^T\bv^k +b^k\pi^T{\bf 1}, \quad \forall k\ge0\ \as
\]
 Since $\sum_{k=0}^\infty a^k= \infty$, and $\sum_{k=0}^\infty b^k<\infty$ {\it a.s.},
the conditions of Lemma~\ref{Lemma-polyak_2} are satisfied
with $v^k = \pi^T\bv^k$, $\a^k=0$, $q^k=a^k$, and $p^k = b^k \pi^T{\bf 1}$, implying \as
$\lim_{k\to\infty}\pi^T\bv^k=0$.  $\{\bv^k\}$ being nonnegative and $\pi>0$ imply
  $\lim_{k\to\infty}\bv^k=0$ \as
\end{proof}

We now proceed to analyze the convergence of Algorithm~2. Defining
 $( \zeta_w^k)^T=\left[(\zeta_{w1}^k)^T,\cdots,(\zeta_{wm}^k)^T\right]$ with $
 \zeta_{wi}\triangleq \sum_{j\in\mathbb{N}_{R,i}^{\rm in}}R_{ij}\zeta_j^k
$ and  $( \xi_w^k)^T=\left[(\xi_{w1}^k)^T,\cdots,(\xi_{wm}^k)^T\right]$ with $
 \xi_{wi}\triangleq \sum_{j\in\mathbb{N}_{C,i}^{\rm in}}C_{ij}\xi_j^k
$, we write the dynamics of Algorithm~2  in the following more compact form:
\begin{equation}\label{eq:push-pull}
\begin{aligned}
x^{k+1}&=\left((I+\gamma_1^kR)\otimes I_d\right) x^k+\gamma_1^k \zeta_w^k-\lambda^k y^k\\
y^{k+1}&=\left(((1-\alpha^k)I+\gamma_2^kC)\otimes I_d\right) y^k +\gamma_2^k \xi_w^k +g^{k+1} \\
&\qquad-(1-\alpha^k)g^{k}
\end{aligned}
\end{equation}
where we used $g^{k+1}=\nabla f(x^{k+1})$ for notational simplicity.

\begin{Lemma 1}\label{lemma:left_right_eigenvectors} \cite{horn2012matrix} (or Lemma 1 in \cite{pu2020push})
Under Assumption \ref{Assumption:push_pull topology},
for every $k$, the matrix $I+\gamma_1^kR$ has a unique
  nonnegative left eigenvector $u^T$ (associated with eigenvalue 1) satisfying $u^T{\bf 1}=m$, and   the matrix $(1-\alpha^k)I+\gamma_2^kC$ has a unique   nonnegative right eigenvector $v$ (associated with eigenvalue $1-\alpha^k$) satisfying $ {\bf 1}^Tv=m$.
\end{Lemma 1}

According to Lemma 3 in~\cite{pu2020push}, we know that the spectral radius of $R^k\triangleq I+\gamma_1^k R-\frac{ {\bf 1}u^T}{m}$ is equal to $1-\gamma_1^k|\nu_R|<1$ where $\nu_R$ is an eigenvalue of $R$. Furthermore,   there exists a vector norm $\|x\|_R\triangleq\|\tilde{R}x\|_2$
(where $\tilde{R}$ is determined by $R$~\cite{pu2020push}) such that $\|R^k\|_R<1$
is arbitrarily close to the spectral radius of $R^k$, i.e., $1-\gamma_1^k|\nu_R|<1$. Without loss of generality, we represent this norm as  $\|R^k\|_R=1-\gamma_1^k \rho_R <1$. Similarly,   we have that the spectral radius of $C^k\triangleq (1-\alpha^k)I+\gamma_2^k C-\frac{ v{\bf 1}^T}{m}$ is equal to $1-\alpha^k-\gamma_2^k|\nu_C|<1$ where $\nu_C$ is an eigenvalue of $C$. Furthermore,   there exists a vector norm $\|x\|_C\triangleq\|\tilde{C}x\|_2$   (where $\tilde{C}$ is  determined by $C$~\cite{pu2020push}) such that $\|C^k\|_C<1$
is arbitrarily close to the spectral radius of $C^k$, i.e., $1-\alpha^k-\gamma_2^k|\nu_C|<1$. Without loss of generality, we bound this norm as  $\|C^k\|_C\leq 1-\gamma_2^k \rho_C <1$.

 Defining $\bar{x}^k=\frac{1}{m}\sum_{i=1}^{m} u_ix_i^k$ and $\bar{y}^k=\frac{1}{m}\sum_{i=1}^{m} y_i^k$,
we have
\begin{equation}\label{eq:bar_x_k+1}
\begin{aligned}
  \bar{x}^{k+1}&=\bar{x}^k+\gamma_1^k \bar{\zeta}_w^k-\lambda^k \frac{(u\otimes { I_d})^T}{m}y^k\\
\bar{y}^{k+1}&=(1-\alpha^k)\bar{y}^k +\gamma_2^k \bar{\xi}_w^k+\bar{g}^{k+1}-(1-\alpha^k)\bar{g}^k
\end{aligned}
\end{equation}
with $
\bar{\zeta}_{w}^k= \frac{1}{m}\sum_{i=1}^m  u_i\zeta_{wi}^k$, $\bar{\xi}_w^k=\frac{1}{m}\sum_{i=1}^m\xi_{wi}^k$,  and $ \bar{g}^k=\frac{1}{m}\sum_{i=1}^m g_i^k$.

From (\ref{eq:bar_x_k+1}), we can further obtain
\begin{equation}\label{eq:x_bar_update}
\begin{aligned}
 \bar{x}^{k+1}&
 =\bar{x}^k+\gamma_1^k \bar{\zeta}_w^k-\lambda^k \frac{(u\otimes { I_d})^T}{m}\left(y^k-(v\otimes {I_d})\bar{y}^k\right)\\
 &\qquad-\lambda^k \frac{(u\otimes {I_d})^T}{m}(v\otimes {I_d})\bar{y}^k
\end{aligned}
\end{equation}
Using the relationship $\lambda^k \frac{(u\otimes {I_d})^T}{m}(v\otimes {I_d})\bar{y}^k=\lambda^k\frac{u^Tv}{m}\bar{y}^k$, we can rewrite (\ref{eq:x_bar_update}) as follows
\begin{equation}\label{eq:x_bar_update2}
\bar{x}^{k+1}\hspace{-0.06cm}= \hspace{-0.06cm} \bar{x}^k-\lambda^k \hspace{-0.06cm}\frac{(u\otimes {I_d })^T}{m}\hspace{-0.1cm}\left(y^k-(v\otimes {I_d})\bar{y}^k\right)-\lambda^k\frac{u^Tv}{m}\bar{y}^k+\gamma_1^k \bar{\zeta}_w^k
\end{equation}

 In what follows, we use $F^*$ to denote the optimal value of the problem in~\eqref{eq:optimization_formulation1},
i.e.,
$F^*=\min_{\theta\in\R^d}F(\theta)$.

Next, we provide a generic convergence result for dynamic-consensus (gradient-tracking) based distributed algorithms for problem~(\ref{eq:optimization_formulation1}). To this end, we need a  measure under the $\|\cdot\|_R$ norm for the distance between  all $x_1^k,\,x_2^k,\cdots,x_m^k$  and $\bar{x}^k$. Following \cite{pu2020push}, we   define a matrix norm for all $x$ iterates $ {\bf x}^k\triangleq \left[x_1^k,\,x_2^k,\cdots,x_m^k\right]^T\in\mathbb{R}^{m\times d}$:
\begin{equation}\label{eq:matrix_norm_R}
\|{\bf x}^k\|_R=\left\|\left[\|{\bf x}^k_{(1)}\|_R,\,\|{\bf x}^k_{(2)}\|_R,\cdots,\|{\bf x}^k_{(d)}\|_R \right] \right\|_2
\end{equation}
 where the subscript $2$ denotes the $2-$norm and ${\bf x}^k_{(i)}$ denotes the $i$th column of ${\bf x}^k$. Defining ${\bf\bar{x}}^k$ as $\left[\bar{x}^k,\,\bar{x}^k,\cdots,\bar{x}^k\right]^T\in\mathbb{R}^{m\times d}$, one can easily see that $\|{\bf x}^k-{\bf\bar{x}}^k \|_R$ measures the distance between all $x_i^k$ and their average $\bar{x}^k$. Similarly, we  define a matrix norm $\|\cdot\|_C$ for $ {\bf y}^k\triangleq \left[y_1^k,\,y_2^k,\cdots,y_m^k\right]^T\in\mathbb{R}^{m\times d}$:
\begin{equation}\label{eq:matrix_norm_C}
\|{\bf y}^k\|_C=\left\|\left[\|{\bf y}^k_{(1)}\|_C,\,\|{\bf y}^k_{(2)}\|_C,\cdots,\|{\bf y}^k_{(d)}\|_C \right] \right\|_2
\end{equation}
and use $\|{\bf y}^k-{\rm diag}(v){\bf\bar{y}}^k \|_C$ (with ${\rm diag}(v)={\rm diag}(v_1,\, \ldots,v_m)$ and ${\bf\bar{y}}^k\triangleq\left[\bar{y}^k,\, \ldots,\bar{y}^k\right]^T\in\mathbb{R}^{m\times d}$) to measure  the distance between all $y$ iterates and their $v$-weighted average  $\bar{y}^k$.

\begin{Lemma 1}\label{Theorem:general_gradient_tracking}
Assume that the objective function $F(\cdot)$ is
differentiable and that the problem~(\ref{eq:optimization_formulation1}) has an optimal solution.
Suppose that a distributed algorithm generates
sequences $\{x_i^k\}\subseteq\mathbb{R}^d$ and
$\{y_i^k\}\subseteq\mathbb{R}^d$ under coupling matrices $R$ and $C$, respectively,   such that the following relation holds
\as for some sufficiently large integer $T\ge0$ and for all $k\ge T$:
\begin{equation}\label{eq-fine}
\mathbb{E}\left[\bv^{k+1}|\mathcal{F}^k\right] \le \left(V^k\hspace{-0.07cm} + a^k {\bf 1}{\bf 1}^T\right)\bv^{k}+b^k{\bf
1} - H^k \hspace{-0.07cm}\left[\begin{array}{c}
 \|\nabla F(\bar x^k)\|^2\cr
 \|\bar {y}^k\|^2\end{array}\right]
\end{equation}
 where $\mathcal{F}^k=\{x_i^\ell,y_i^\ell;0\le \ell\le k,\ i\in[m]\}$   and
\[
\begin{aligned}
&\bv^k=\left[\begin{array}{c}\bv^k_1\\\bv^k_2\\\bv^k_3\end{array}\right]\triangleq\left[\begin{array}{c}  F(\bar
x^{k})-F^* \cr  \|{\bf x}^{k}-{\bf \bar x}^{k}\|_R^2\cr
\|{\bf y}^{k}-{\rm diag}(v){\bf \bar y}^{k}\|_C^2\end{array}
 \right],\\
&\hspace{-0.2cm}
 \begin{array}{ll}
 V^k=& H^k=\\
  \left[\begin{array}{ccc} 1 &  \kappa_1\lambda^k &\kappa_2\lambda^k\\
                            0 & 1-\kappa_3\gamma_1^k&0\\
                            0 & 0&1-\kappa_4\gamma_2^k
 \end{array}\right],& \hspace{-0.1cm}
 \left[\begin{array}{cc} \kappa_5\lambda^k & \kappa_6 \lambda^k-\kappa_7(\lambda^k)^2\\0 &
-\kappa_8\frac{(\lambda^k)^2}{\gamma_1^k}\cr 0 & -\kappa_9\frac{(\lambda^k)^2}{\gamma_2^k}
\end{array} \right]
 \end{array}
 %
\end{aligned}
\] with  $\kappa_i>0$ for all  $1\leq i\leq 9$  and $\kappa_3,\kappa_4\in(0,1)$, while the nonnegative scalar
sequences $\{a^k\}$, $\{b^k\}$ and positive sequences  $\{\lambda^k\}$,    $\{\gamma_1^k\}$, $\{\gamma_2^k\}$ satisfy
$\sum_{k=0}^\infty a^k<\infty$ {\it a.s.}, $\sum_{k=0}^\infty b^k<\infty$ {\it a.s.}, $\sum_{k=0}^\infty \lambda^k=\infty$,  $\sum_{k=0}^\infty \gamma_{\imath}^k=\infty$, $\sum_{k=0}^\infty (\gamma^k_\imath)^2<\infty$,   $\sum_{k=0}^\infty \frac{(\lambda^k)^2}{\gamma_\imath^k}<\infty$, $\lim_{k\to\infty}\lambda^k/\gamma_\imath^k=0$ for $\imath\in\{1,\,2\}$, and $\lim_{k\to\infty}\gamma_1^k/\gamma_2^k<\infty$.
Then, we have:
\begin{itemize}
\item[(a)]  $\lim_{k\to\infty} F(\bar x^k)$
exists \as and
\[
\lim_{k\to\infty}\|x_i^k - \bar x^k\|= \lim_{k\to\infty}\|y_i^k -
v_i\bar y^k\|=0, \: \forall i\quad a.s. \]
\item[(b)]
$\lim\inf_{k\rightarrow\infty}\|\nabla F(\bar x^k)\|=0$ holds \as
 Moreover, if the function $F(\cdot)$ has bounded level sets, then
$\{\bar x^k\}$ is bounded and every accumulation point of  $\{\bar
x^k\}$ is an optimal solution {\it a.s.},
and
$\lim_{k\rightarrow\infty}F(x_i^k)= F^*$ \as for all $i\in[m]$.
\end{itemize}

\end{Lemma 1}
\begin{proof}
See Appendix C.
\end{proof}

\begin{Remark 1}
In Lemma~\ref{Theorem:general_gradient_tracking}(b),  the bounded level set condition
can be replaced with
any other condition ensuring that the sequence $\{\bar x^k\}$ is \as bounded.
\end{Remark 1}



 Lemma~\ref{Theorem:general_gradient_tracking} is critical for establishing
convergence properties of the gradient tracking-based distributed algorithm together
with suitable conditions on the DP-noise injected by the agents. We make the following assumption on the noise:
\begin{Assumption 1}\label{assumption:dp-noises-intrack}
For every $i\in[m]$, the noise sequences $\{\zeta_i^k\}$ and
$\{\xi_i^k\}$ are zero-mean independent random variables, and independent of $\{x_i^0;i\in[m]\}$.
Also, for every $k$, the noise collection $\{\zeta_j^k,\xi_j^k; j\in[m]\}$ is independent. The noise variances
$(\sigma_{\zeta,i}^k)^2=\mathbb{E}\left[\|\zeta_i^k\|^2\right]$ and $(\sigma_{\xi,i}^k)^2=\mathbb{E}\left[\|\xi_i^k\|^2\right]$
and their attenuation stepsizes $\g_1^k$ and $\g_2^k$ are
such that
\begin{equation}\label{eq:condition_assumption5}
\begin{aligned}
\sum_{k=0}^\infty(\g_1^k)^2&\max_{i\in[m]}(\sigma_{\zeta,i}^k)^2<\infty,  \sum_{k=0}^\infty (\gamma_2^k)^2\, \max_{j\in[m]}(\sigma_{\xi,j}^k)^2 <\infty.
\end{aligned}
\end{equation}
The initial random vectors satisfy
$\mathbb{E}\left[\|x_i^0\|^2\right]<\infty$,  $\forall i\in[m]$.
\end{Assumption 1}
{
\begin{Remark 1}
Given that $\gamma_1^k$, $\gamma_2^k$, and $\lambda^k$ decrease  with time, (\ref{eq:condition_assumption5}) can be satisfied even when $\{\sigma_i^k\}$ increases with time. For example, under $\lambda^k=\mathcal{O}(\frac{1}{k})$, $\gamma_1^k=\mathcal{O}(\frac{1}{k^{0.9}})$, $\gamma_2^k=\mathcal{O}(\frac{1}{k^{0.7}})$, an increasing $\{\sigma_i^k\}$ with increasing rate no faster than $\mathcal{O}(k^{0.15})$
still satisfies the summable condition in (\ref{eq:condition_assumption5}).
%
\end{Remark 1}
}
\begin{Assumption 1}\label{as:bounded_gradients}
 The gradients of all individual objective functions are bounded, i.e., there exists a constant $C$ such that $\|\nabla f_i(x)\|_1\leq C$ holds for all $x\in\mathbb{R}^p$ and $1\leq i \leq m$.
\end{Assumption 1}

\begin{Theorem 1}\label{theorem_convergence_algorithm2}
 Let  Assumption 1, Assumption~\ref{Assumption:push_pull topology}, Assumption~\ref{assumption:dp-noises-intrack}, and Assumption \ref{as:bounded_gradients}
hold.
If nonnegative sequences $\{\gamma_1^k\}$, $\{\gamma_2^k\}$, $\{\alpha^k\}$, and $\{\lambda^k\}$ satisfy
$
 \sum_{k=0}^\infty \gamma_\imath^k=\infty$,  $\sum_{k=0}^\infty (\gamma_\imath^k)^2<\infty$, $\sum_{k=0}^\infty \alpha^k=\infty$,  $\sum_{k=0}^\infty \lambda^k=\infty$, $\sum_{k=0}^\infty \frac{(\lambda^k)^2}{\gamma_\imath^k}<\infty$, $\lim_{k\to\infty}\lambda^k/\gamma_\imath^k=0$ for $\imath\in\{1,\,2\}$, $\lim_{k\to\infty}\lambda^k/\alpha^k<\infty$, $\sum_{k=0}^\infty \frac{(\alpha^k)^2}{\gamma_2^k}<\infty$ and $\sum_{k=0}^\infty \frac{(\gamma_1^k)^2}{\gamma_2^k}<\infty$,
then, the results of Lemma~\ref{Theorem:general_gradient_tracking}
hold for Algorithm~2.
\end{Theorem 1}
\begin{proof}
See Appendix D.
\end{proof}
\begin{Remark 1}
In networked systems, usually  communication imperfections can be modeled as channel noises \cite{kar2008distributed}, which can be regarded as a special case of the DP noise considered here.  Therefore, Algorithm 2 can also be used to counteract such communication imperfections in distributed optimization.
\end{Remark 1}

{

\begin{Remark 1}
 Because the evolution of $x_i^k$ to the optimal solution satisfies the conditions in Lemma \ref{Theorem:general_gradient_tracking}, which are in turn derived based on Lemma \ref{Lemma-polyak_2}, we can leverage Lemma \ref{Theorem:general_gradient_tracking} and Lemma \ref{Lemma-polyak_2} to characterize the convergence speed.  More specifically, in the proof of Lemma \ref{Theorem:general_gradient_tracking} in the appendix, (\ref{eq-finer}) and the relationship $\tilde \pi^T
\tilde V^k=(1-\a\gamma_1^k) \tilde\pi^T$ imply that $\bv^k_2\triangleq \|{\bf x}^{k}-{\bf \bar x}^{k}\|_R^2$ and $\bv^k_3 \triangleq \|{\bf y}^{k}-{\rm diag}(v){\bf \bar y}^{k}\|_C^2$  decay  to zero with a rate no slower than $\mathcal{O}(\frac{\gamma_1^k}{\gamma_2^k})$. Furthermore, (\ref{eq-sumfinite}) implies that $\lambda^k\|\nabla F(\bar{x}^k)\|^2$ decays to zero with a rate no slower than $\mathcal{O}(\frac{1}{k})$, i.e., $\|\nabla F(\bar{x}^k)\|^2$ decays to zero with a rate no slower than $\mathcal{O}(\frac{1}{k\lambda^k})$.   Moreover, from the proof in Lemma \ref{Theorem:general_gradient_tracking} (specifically (\ref{eq-finer}) and the paragraph below it), we know that the decreasing speed of $\|{\bf x}^{k}-{\bf \bar x}^{k}\|_R^2$, $\|{\bf y}^{k}-{\rm diag}(v){\bf \bar y}^{k}\|_C^2$ increases with an increase in $\alpha$, which in turn increases with an increase in $\kappa_3$ and $\kappa_4$. Further noting that $\kappa_3$ and $\kappa_4$  correspond  to the spectral radius of $R$ and $C$, respectively, we have that the convergence speed increases with an increase in the spectral radius of  $R$ and $C$ defined in Assumption \ref{Assumption:push_pull topology} (see  (\ref{eq:diagonal_entries}) for diagonal entries).

\end{Remark 1}

\subsection{Privacy analysis}


\begin{Theorem 1}\label{th:DP_Algorithm2}
Under Assumptions \ref{assumption:objective_functions}, \ref{Assumption:push_pull topology} and \ref{as:bounded_gradients},  if $F(\cdot)$ has bounded level sets, nonnegative sequences $\{\lambda^k\}$, $\{\alpha^k\}$, $\{\gamma_1^k\}$, and $\{\gamma_2^k\}$ satisfy the conditions in Theorem \ref{theorem_convergence_algorithm2}, and all elements of $\zeta_i^k$ and $\xi_i^k$ are drawn independently from Laplace distribution ${\rm Lap}(\nu^k)$ with $(\sigma_{\zeta,i}^k)^2=(\sigma_{\xi,i}^k)^2=2(\nu^k)^2$ satisfying Assumption \ref{assumption:dp-noises-intrack}, then all agents will converge \as to an optimal solution. Moreover,
{
\begin{enumerate}
\item For any finite number of iterations $T$, Algorithm 1 is  $\epsilon$-differentially private with the cumulative privacy budget bounded by $\epsilon\leq \sum_{k=1}^{T}\frac{2C(\varsigma_x^k+\varsigma_y^k)}{\nu^k}$   where $\varsigma_x^k\triangleq \sum_{p=1}^{k-1}(\Pi_{q=p}^{k-1}(1-\bar{R}\gamma_1^{q})) \lambda^{p-1}\varsigma_y^{p-1}+ \lambda^{k-1}\varsigma_y^{k-1}$, $\varsigma_y^k\triangleq \sum_{p=1}^{k-1}(\Pi_{q=p}^{k-1}(1-\alpha^q-\bar{C}\gamma_2^{q}))(2-\alpha^{p-1})+2-\alpha^{k-1}$, $\bar{R}\triangleq\min_i\{|R_{ii}|\}$, $\bar{C}\triangleq\min_i\{|C_{ii}|\}$, and $C$ is from Assumption \ref{as:bounded_gradients};
\item  The cumulative privacy budget is  finite for $T\rightarrow\infty$  when the sequence  $\{\frac{\lambda^k}{\nu^k}\}$ is summable.
\end{enumerate}
}
%
%
\end{Theorem 1}
\begin{proof}
Since the convergence follows Theorem \ref{theorem_convergence_algorithm2}, we only consider the privacy statements.

Following   Definition \ref{de:sensitivity} and the argument in the proof of Theorem \ref{th:DP_Algorithm1}, we know that the sensitivity of Algorithm 2 is determined by $x_i^k-{x_i'}^k$ and $y_i^k-{y_i'}^k$, which, according to Algorithm 2,  have the following dynamics:
 \[
 \begin{aligned}
x_i^{k+1}-{x'_i}^{k+1}=&(1-\gamma_1^k |R_{ii}|)(x_i^k-{x'_i}^k)-\lambda^k(y_i^k-{y'_i}^k),\\
y_i^{k+1}-{y'_i}^{k+1}=&(1-\alpha^k-\gamma^k_2|C_{ii}|)(y_i^k-{y'_i}^k)\\
 & +(g_i^{k+1}-{g'_i}^{k+1})- (1-\alpha^k)(g_i^k-{g'_i}^k),
 \end{aligned}
 \]
 where we have represented $\nabla f_i(x_i^k)$ and $\nabla f'_i({x'_i}^k)$ as $g_i^k$ and ${g'_i}^k$, respectively,  for notational simplicity.
 Note that we have also used  the fact that the observations $x_j^k+\zeta_j^k$ (resp. $y_j^k+\xi_j^k$) and ${x'_j}^k+{\zeta'_j}^k$ (resp. ${y'_j}^k+{\xi'_j}^k$) are the same.

 Hence,   $y_i^k-{y'_i}^k$ satisfies
 \[
 \begin{aligned}
 &\|y_i^{k+1}-{y'_i}^{k+1}\|_1\\
 &\leq (1-\alpha^k-\gamma^k_2|C_{ii}|)\|y_i^k-{y'_i}^k\|_1\\
 &\qquad + \|g_i^{k+1}-{g'_i}^{k+1}\|_1+(1-\alpha^k)\|g_i^k-{g'_i}^k\|_1\\
 & \leq (1-\alpha^k-\gamma^k_2|C_{ii}|)\|y_i^k-{y'_i}^k\|_1 + (2-\alpha^k)2C,
 \end{aligned}
 \]
 where $C$ is from Assumption \ref{as:bounded_gradients}.

By iteration, we have that  $\|y_i^k-{y'_i}^k\|_1$ is always bounded by $2C\varsigma_y^k$ in the first privacy statement.

 One can also see that  $x_i^k-{x'_i}^k$ satisfies
  \[
 \begin{aligned}
 \|x_i^{k+1}-{x'_i}^{k+1}\|_1&\leq (1 -\gamma^k_1|R_{ii}|)\|x_i^k-{x'_i}^k\|_1\\
 & + \lambda^k\|y_i^k-{y'_i}^k\|_1.
 \end{aligned}
 \]
 Hence, using iteration, we obtain that $\|x_i^k-{x'_i}^k\|_1$ is always bounded by $2C\varsigma_x^k$ in the first privacy statement.

 Therefore, the sensitivity  at iteration $k$ is no larger than $2C(\varsigma_x^k+\varsigma_y^k)$ in the first privacy statement, and, hence, we have the first privacy statement on the cumulative privacy budget for a finite number of iterations.

 On the infinite time horizon, we follow  the argument in the proof of Theorem \ref{th:DP_Algorithm1}. More specifically,  we can prove that the sequence $\{\|g_i^k-{g'_i}^k\|_1\}$ can be bounded by the sequence $\{C_g\gamma_2^k\lambda^k\}$ (with $C_g$ some constant) using the third condition in Definition \ref{de:adjacency} and the guaranteed convergence. Hence, according to Lemma \ref{le:chung},   there always exists a constant $\bar{C}_y$ such that  $ \|y_i^k-{y'_i}^k\|_1\leq \bar{C}_y\lambda^k$ holds. Still using Lemma \ref{le:chung}, we can prove that there always exists a constant $\bar{C'_x}$ such that  $ \|x_i^k-{x'_i}^k\|_1\leq \bar{C'_x}\frac{(\lambda^k)^2}{\gamma_1^k}$ holds. Given that $\lambda^k$ decreases faster than $\gamma_1^k$, we have $ \|x_i^k-{x'_i}^k\|_1\leq \bar{C}_x  \lambda^k  $ for some constant $\bar{C}_x$. Therefore, on the infinite time horizon, the sensitivity is on the order of $\lambda^k$. Hence, we have the result on the cumulative privacy budget when $T\rightarrow\infty$ in the second statement.
\end{proof}
\begin{Remark 1}
 Since we use the standard $\epsilon$-DP framework, we characterize the cumulative privacy budget directly.  Under  relaxed (approximate) $\epsilon$-DP frameworks, such as $(\epsilon,\delta)$-DP \cite{kairouz2015composition}, zero-concentrated DP \cite{bun2016concentrated}, or R\'{e}nyi DP \cite{mironov2017renyi},    advanced composition theories in \cite{kairouz2015composition,bun2016concentrated,mironov2017renyi}  can be exploited to characterize the cumulative privacy budget.
\end{Remark 1}
}
\section{Numerical Experiments}\label{sec-numerics}
{
\subsection{Evaluation using distributed estimation}
We first evaluate  the performance of the two proposed    algorithms using a canonical distributed estimation problem}
where a
 network of $m$ sensors collectively estimate an  unknown parameter
$\theta\in\mathbb{R}^d$. More specifically, we assume that each
sensor $i$ has a noisy measurement  of the parameter,
$z_{i}=M_i\theta+w_{i}$, where $M_i\in\mathbb{R}^{s\times d}$ is the
measurement matrix of agent $i$ and $w_{i}$ is Gaussian measurement noise of unit variance.
Then the maximum likelihood  estimation of parameter $\theta$  can be solved using the
optimization problem formulated as
(\ref{eq:optimization_formulation1}), with each $f_i(\theta)$ given
as
$
f_i(\theta)=\|z_i-M_i\theta\|^2+\varsigma\|\theta\|^2
$
where $\varsigma$ is a  regularization parameter \cite{xu2017convergence}.

\begin{figure}
\center
\includegraphics[width=0.22\textwidth]{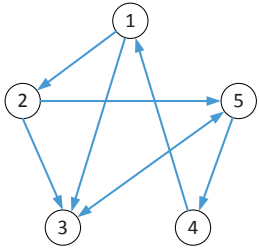}
    \caption{The interaction topology of the network}
    \label{fig:topology}
\end{figure}

\begin{figure}
\includegraphics[width=0.5\textwidth]{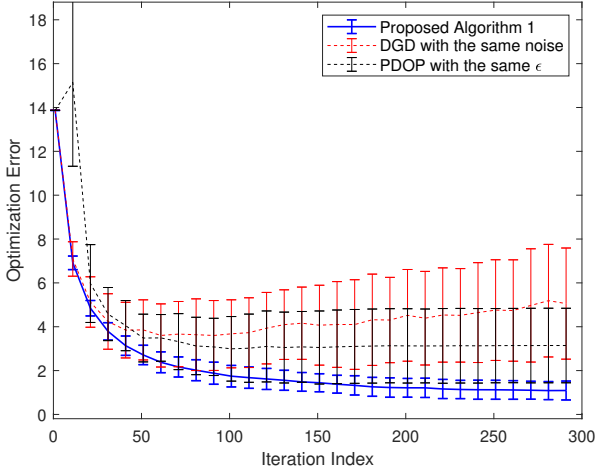}
    \caption{Comparison of Algorithm 1 with existing  distributed gradient descent algorithm (DGD) in \cite{nedic2009distributed} (under the same noise) and the differential-privacy approach for decentralized optimization PDOP in \cite{huang2015differentially} (under the same privacy budget) using the distributed estimation problem}
    \label{fig:comparison_DGD}
\end{figure}

We consider a network of $m=5$ sensors interacting on the graph depicted in Fig. \ref{fig:topology}. In the evaluation, we set $s=3$ and $d=2$. To evaluate the performance of the proposed Algorithm 1, we ignored the directions of edges in Fig.  \ref{fig:topology} in the selection of coupling weights and injected Laplace based DP-noise with parameter $\nu^k=1+0.1k^{0.3}$  in every message shared in all iterations. We set the stepsize $\lambda^k$ and diminishing sequence $\gamma^k$ as $\lambda^k=\frac{0.02}{1+0.1k}$ and $\gamma^k=\frac{1}{1+0.1k^{0.9}}$, respectively, which satisfy the conditions in Theorem 1 and Theorem 2. In the evaluation, we ran our algorithm for 100 times and calculated the average as well as the variance of the optimization error as a function of the iteration index. The result is given by the blue curve and error bars in Fig. \ref{fig:comparison_DGD}. For comparison, we also ran the existing static-consensus based distributed gradient descent (DGD) approach in \cite{nedic2009distributed} under the same noise, and the differential-privacy approach for distributed optimization  (PDOP) in \cite{huang2015differentially} under the same privacy budget. Note that PDOP uses  geometrically decreasing stepsizes (which are summable) to ensure a finite privacy budget, but the fast decreasing stepsize also leads to optimization errors.   The evolution of the average optimization error and variance of the DGD and PDOP approaches are given by  the red and black curves/error bars in Fig. \ref{fig:comparison_DGD}, respectively. It is clear that the proposed algorithm has a comparable convergence speed but much better optimization accuracy.

\begin{figure}
\includegraphics[width=0.5\textwidth]{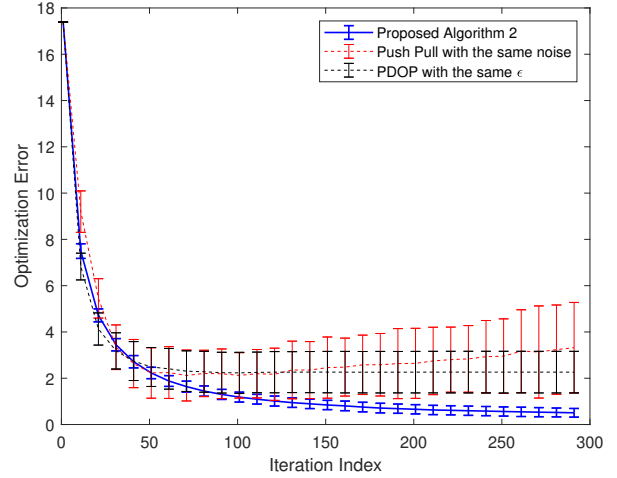}
    \caption{Comparison of Algorithm 2 with existing dynamic-consensus based distributed gradient algorithm (Push Pull) in \cite{pu2020push} (under the same noise) and the PDOP-based differential-privacy approach in \cite{huang2015differentially} for Push Pull (under the same privacy budget) using the distributed estimation problem}
    \label{fig:comparison_Push_pull}
\end{figure}

We also evaluated Algorithm 2 which is applicable to general directed graphs. More specifically, still using the topology in Fig. \ref{fig:topology}, we selected $R$ and $C$ matrices according to  Assumption \ref{Assumption:push_pull topology} and  set the stepsize and diminishing sequences as $\lambda^k=\frac{0.02}{1+0.1k}$, $\alpha^k=\frac{0.02}{1+0.1k}$ $\gamma_1^k=\frac{1}{1+0.1k^{0.9}}$, and $\gamma_2^k=\frac{1}{1+0.1k^{0.7}}$, respectively. We injected Laplace noises $\zeta_i^k$ and $\xi_i^k$ (both have  parameter $\nu^k= 1+0.1k^{0.1}$) on all shared $x_i^k$ and $y_i^k$ respectively to enable DP, and it can be verified that the parameters satisfy the conditions in Theorem 3 and Theorem 4. We ran our algorithm for 100 times and calculated the average as well as the variance of the optimization error as a function of the iteration index.  The result is given by the blue curve and error bars in Fig. \ref{fig:comparison_Push_pull}. For comparison, we also ran the conventional dynamic-consensus based Push-Pull method in \cite{pu2020push} under the same noise and the PDOP based differential-privacy approach for distributed optimization. Because the PDOP based approach requires the stepsize to  decay  with a geometric rate, we set the stespize of Push-Pull to $0.95^k$ and used a geometrically decaying noise such that it has the same privacy budget as our approach. The evolution of the average optimization error and variance of Push Pull (with the same noise as our approach) and PDOP-privacy based Push Pull (with  the same privacy budget as our approach) are depicted by  the red and black curves/error bars in Fig. \ref{fig:comparison_Push_pull}, respectively. It is clear that the proposed algorithm has a comparable convergence speed but gained significant improvement in optimization accuracy.

\begin{figure}
\includegraphics[width=0.5\textwidth]{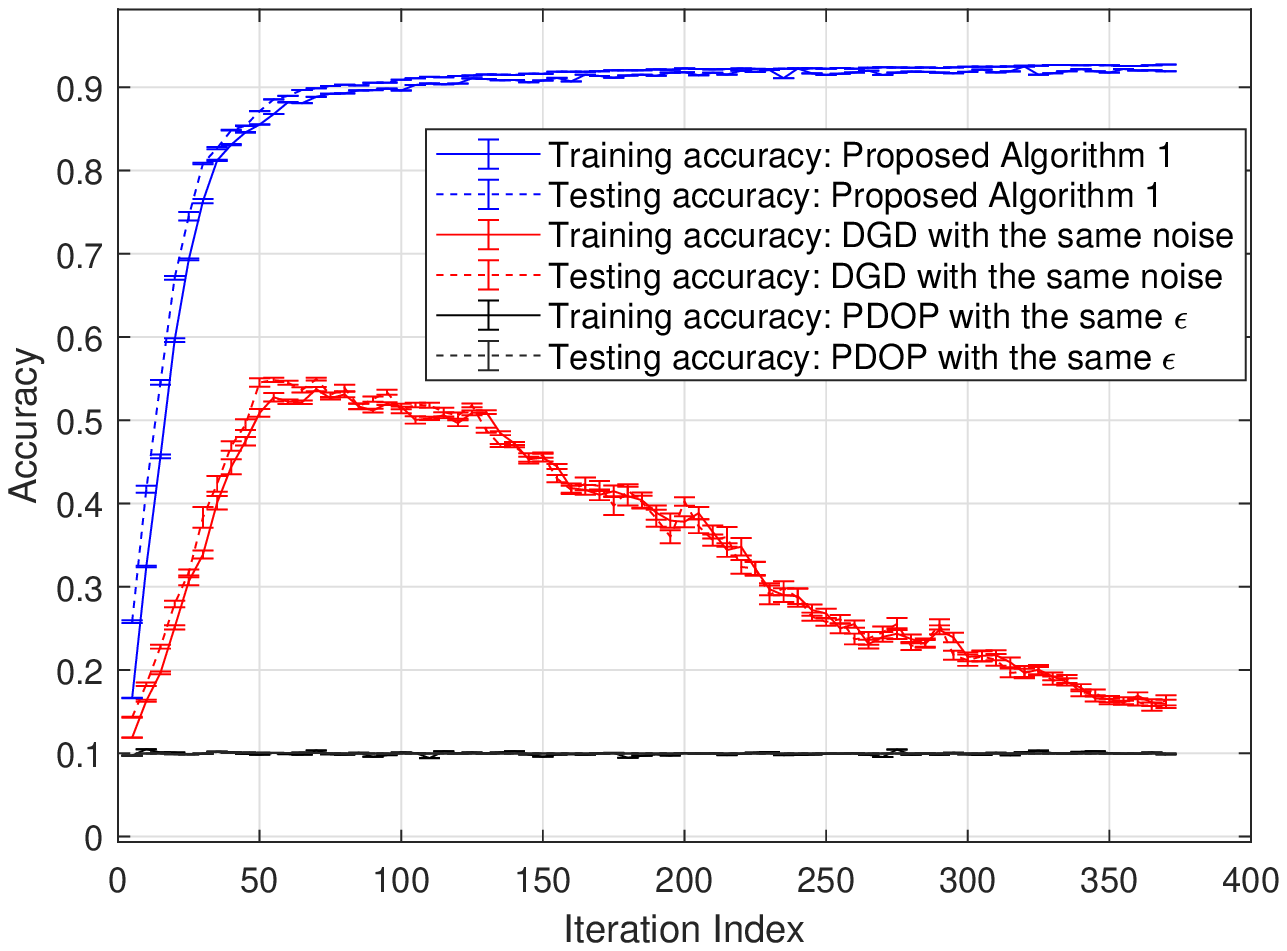}
    \caption{Comparison of Algorithm 1 with existing  distributed gradient descent algorithm (DGD) in \cite{nedic2009distributed} (under the same noise) and the differential-privacy approach for decentralized optimization PDOP in \cite{huang2015differentially} (under the same privacy budget) using the MNIST image classification problem}
    \label{fig:comparison_MNIST1}
\end{figure}

\begin{figure}
\includegraphics[width=0.5\textwidth]{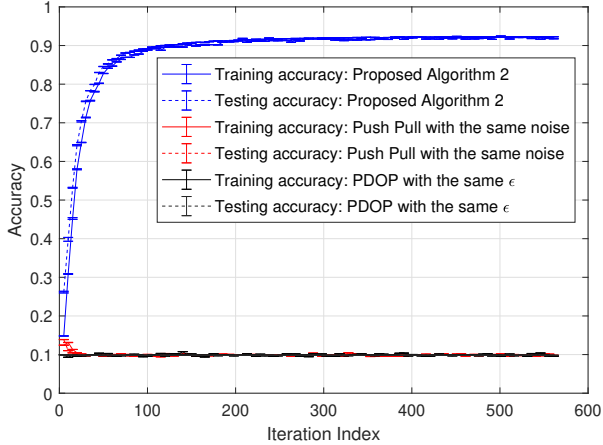}
    \caption{Comparison of Algorithm 2 with existing dynamic-consensus based distributed gradient algorithm (Push Pull) in \cite{pu2020push} (under the same noise) and the PDOP-based differential-privacy approach in \cite{huang2015differentially} for Push Pull (under the same privacy budget) using the MNIST image classification problem}
    \label{fig:comparison_MNIST2}
\end{figure}
{
\subsection{Evaluation using image classification on MNIST}
We also used decentralized training of a convolutional neural network
(CNN)  to evaluate the performance of our proposed algorithms.
More specially, we   consider five agents which collaboratively
train a CNN using the MNIST dataset \cite{MNIST} under the topology
 in Fig. \ref{fig:topology}. The MNIST  data set is a large
benchmark database of handwritten digits widely used for training
and testing in the field of machine learning \cite{deng2012mnist}.
Each agent has a local copy of the CNN. The CNN has 2 convolutional
layers with 32 filters with each followed by a max pooling layer,
and then two more convolutional layers with 64 filters each followed
by another max pooling layer and a dense layer with 512 units. Each
agent has access to a portion of the MNIST dataset, which was
further divided into two subsets for training and validation,
respectively. To evaluate the proposed Algorithm 1, We set the stepsize as  $\lambda^k= \frac{1}{1+0.01k}$ and the weakening factor $\gamma^k$ as $\frac{1}{1+0.01k^{0.9}}$. The Laplace noise parameter was set to $\nu^k=1+0.01k^{0.3}$ to enable $\epsilon$-DP.
The evolution of the training and testing accuracies averaged
over 50 runs are illustrated by the solid and dashed blue curves
in Fig. \ref{fig:comparison_MNIST1}. To compare the convergence performance
of our algorithm with the conventional distributed
gradient descent algorithm under DP-noise, we also implemented the distributed
 gradient descent (DGD)  algorithm in \cite{nedic2009distributed} to train
the same CNN {using stepsize $\frac{1}{1+0.01k}$} under the same Laplace noise. The results are illustrated by the solid and dotted red curves in Fig. \ref{fig:comparison_MNIST1}. It can be seen that the proposed algorithm has  much better robustness to DP-noise.  Moreover, to compare with the existing DP approach for distributed optimization, we also implemented the DP approach PDOP in \cite{Huang15}  on DGD   under the same privacy budget $\epsilon$. PDOP uses geometrically decaying stepsizes and noises to ensure a finite privacy budget. However,   such fast-decaying stepizes  turned out to be unable to train the complex CNN model (see training and testing accuracies in solid and dashed black curves in Fig. \ref{fig:comparison_MNIST1}, respectively under   $\lambda^k=0.95^k$ and $\nu^k=0.98^k$). These comparisons corroborate the advantage of the proposed Algorithm 1.

To show the influence of DP-noise on the final optimization accuracy, we also scaled the noise by $0.5$ and $2$ respectively and obtained the training and testing accuracies. To compare the strength of enabled privacy protection, we ran the
 DLG attack model proposed in \cite{zhu2019deep}, which is the most
powerful inference algorithm reported to date in terms of reconstructing exact raw data from shared gradient/model updates. The attacker was assumed to be able to observe
all messages shared among the agents. The training/testing accuracies under different levels of DP-noise and DLG attacker's inference errors are summarized in Table 1.  It can be seen that there is a trade-off between privacy and accuracy under a fixed iteration number $20,000$.
\begin{table*}
{
{\caption {Training/testing accuracies and DLG
attacker's inference errors under different levels of
DP-noise}}
 \center
 \noindent\makebox[\textwidth][c]{\hspace{2cm}
\begin{minipage}{13cm}
\begin{tabular}{  c|c|c|c|c|c|c  }
 &\multicolumn{2}{r}{Algorithm 1}  & & \multicolumn{2}{r} {Algorithm 2}&\\
  \hline
 Noise Level\footnote{Considering the noise in Fig. 4 and Fig. 5 as the base level for Algorithm 1 and Algorithm 2, respectively. } &$\times 0.5$ & $\times 1$ & $\times 2$&$\times 0.5$& $\times 1$&$\times 2$ \\
  \hline
  Training Accuracy   &0.951 & 0.925& 0.859& 0.924 & 0.921 &0.910\\
  \hline
  Testing Accuracy &0.951 & 0.929& 0.861& 0.926 & 0.922 & 0.913\\
  \hline
  Final DLG Error & 310.2& 350.3 & 412.5 & 301.1 &336.7&389.7\\
  \hline
\end{tabular}
\end{minipage}}}
\end{table*}

Using the same interaction topology, CNN network, and MNIST dataset, we also evaluated the performance of the proposed Algorithm 2 under DP-noise. The parameters of Algorithm 2 were set as $\lambda^k= \frac{1}{1+0.01k}$, $\alpha^k= \frac{0.1}{1+0.01k}$, $\gamma_1^k= \frac{1}{1+0.01k^{0.9}}$, and $\gamma_2^k= \frac{1}{1+0.01k^{0.7}}$. The Laplace noise parameter was set as  $\nu^k= \frac{1}{1+0.01k^{0.1}}$. The evolution of the training and testing accuracies averaged
over 50 runs are illustrated by the solid and dashed blue curves
in Fig. \ref{fig:comparison_MNIST2}. For comparison, we also implemented the dynamic-consensus based Push Pull algorithm in \cite{pu2020push} to train
the same CNN {using stepsize $0.02$} under the same Laplace noise. The results are illustrated by the solid and dotted red curves in Fig. \ref{fig:comparison_MNIST2}. It can be seen that the same amount of noise, which is tolerable to our proposed Algorithm 2, completely prevents the Push Pull algorithm from training the CNN model.  Moreover, we also applied PDOP based   DP approach in \cite{Huang15} to Push Pull, which  uses geometrically decaying stepsizes and noises to ensure a finite privacy budget.  However, under the same privacy budget, the fast-decaying stepize for Push Pull turned out to be unable to train the complex CNN model either (see Fig. \ref{fig:comparison_MNIST2} for training and testing accuracies in solid and dashed black curves, respectively, under   $\lambda^k=0.95^k$ and $\nu^k=0.98^k$). These comparisons corroborate the advantage of the proposed Algorithm 2.

To   show the influence of DP-noise on the final optimization accuracy and the strength of enabled privacy, we also  scaled the noise by $0.5$ and $2$ respectively and obtained the training/testing accuracies as well as DLG attacker's inference errors. The results are given in Table 1, which shows a trade-off between privacy and accuracy under a fixed iteration number $20,000$.   }

\section{{Conclusions, Discussions, and Future Work}}\label{sec-concl}
Although DP is becoming the de facto standard for publicly sharing information, its direct incorporation into distributed optimization leads to significant reduction in optimization accuracy due to the need to iteratively and repeatedly inject  independent noises.
This paper proposes two DP-oriented gradient based distributed optimization algorithms that ensure both {$\epsilon$-DP} and optimization accuracy. Specifically, the two algorithms can ensure almost sure convergence of all agents to the optimal solution even in the presence of  persistent DP noise. {Both algorithms are also proven able to ensure   $\epsilon$-DP with a finite cumulative privacy budget,  even when the number of iterations goes to infinity. The simultaneous achievement of both provable convergence to the accurate solution and rigorous $\epsilon$-DP with guaranteed finite cumulative privacy budget, to our knowledge, has not been reported before in distributed optimization.}    Numerical simulations {and experimental results using a benchmark dateset} confirm that both algorithms have a better accuracy compared with their respective existing counterparts, while maintaining a comparable  convergence speed.

It is worth noting that our algorithms' simultaneous achievement of both provable convergence to the optimal solution and $\epsilon$-DP does not contradict the fundamental theory and limitations of DP in \cite{dwork2014algorithmic}. Firstly, our convergence guarantee (almost sure convergence) is obtained in the stochastic sense, which is different from deterministic convergence under no DP noise. More specifically, when   the number of implementations tends to infinity, the concept of almost sure convergence   still allows for a finite number of implementations that do not converge to the optimal solution. Secondly, according to the DP theory, conventional query mechanisms on a dataset can achieve  $\epsilon$-DP only by sacrificing query accuracies, but the distributed optimization algorithm does  not correspond to a simple query mechanism on the optimal solution. Instead, what are queried in every iteration of distributed optimization are  individual objective functions (gradients), and revealing the precise optimal solution is not equivalent to revealing accurate objective functions (the actual query target). In fact,  in the language of machine learning, distributed optimization can be viewed as the empirical risk minimization problem,  and the obtained optimal solution corresponds to the optimal model parameter in machine learning. On pages 216-218 of \cite{dwork2014algorithmic}, the authors explicitly  state  that ``the constraint of privacy is not necessarily at odds with the goals of machine learning, both of which aim to extract information from the distribution from which the data was drawn, rather than from individual data points," and ``we are often able to perform private machine learning nearly as accurately, with nearly the same number of examples, as we can perform non-private machine learning." Actually,  under Valiant’s model of machine learning (PAC), \cite{dwork2014algorithmic} notes that a model parameter (called function in \cite{dwork2014algorithmic}) is PAC learnable if and only if it is  PAC learnable under DP (see page 221  of \cite{dwork2014algorithmic}).
Thirdly, the achievement of $\epsilon$-DP does incur utility cost. More specifically, in terms of Algorithm 1, in order to  reduce $\epsilon$ to enhance privacy, we can use a faster-increasing $\{\nu^k\}$ according to Theorem 2, which requires $\{\gamma^k\}$ to decrease faster according to Assumption \ref{assumption:dp-noise}. Given that   the convergence speed is determined by $\mathcal{O}((\frac{\lambda^k}{\gamma^k})^{0.5})$ according to Remark 3, we arrive at the conclusion that a faster decreasing $\{\gamma^k\}$ corresponds to a stronger privacy level but a slower convergence speed. The same conclusion can be drawn for Algorithm 2. In future work, we will systematically quantify the cost of achieving DP  in distributed optimization under the constraint of provable convergence to the optimal solution.  Furthermore, we  also plan to investigate if gradually reducing  communication frequency can   enable rigorous DP.

\section*{Acknowledgement}
The authors would like to thanks Ben Liggett for the help in numerical experiments. They would also like to thank the associate editor and anonymous reviewers, whose comments helped improve  the paper.

\section*{Appendix}
\subsection{Proof of Lemma \ref{th-main_decreasing}}\label{se:Lemma_5}
Let $\theta^*$ be an arbitrary but fixed optimal solution of problem (1). Then, we have
 $F(\bar x^k)- F(\theta^*)\ge0$ for all $k$.
Hence,  by letting $\bv^k=\left[\|\bar x^k-\theta^*\|^2,\ \sum_{i=1}^m \|x_i^k-\bar x^k\|^2\right]^T$,
from relation~(\ref{eq:Theorem_decreasing}) it follows \as that  for all $k\ge0$,
\be\label{eq-fin0}
\hspace{-0.2cm}\mathbb{E}\left[\bv^{k+1}|\mathcal{F}^k\right]
\le \left( \left[\begin{array}{cc}
1 & \frac{\gamma^k}{m}\cr
0& 1-\kappa\gamma^k\cr\end{array}\right] +a^k {\bf 1}{\bf 1}^T \right)\bv^k+b^k{\bf 1}
\ee
Consider the vector $\pi=[1, \frac{1}{m\kappa}]^T$ and note $ \pi^T \left[\begin{array}{cc}
1 & \frac{\gamma^k}{m}\cr
0& 1-\kappa\gamma^k\cr\end{array}\right]= \pi^T$.
Thus, relation~\eqref{eq-fin0} satisfies all   conditions of Lemma~\ref{th-dsystem}.
So it follows that
  $\lim_{k\to\infty}\pi^T\bv^k$ exists {\it a.s.}, and that the sequences $\{\|\bar x^k-\theta^*\|^2\}$
and $\{\sum_{i=1}^m \|x_i^k-\bar x^k\|^2\}$ are bounded {\it a.s.}
From \eqref{eq-fin0} we have the following relation \as  for the second element of $\bv^k$:
\begin{equation}\label{eq:speed_in_lemma5}
\mathbb{E}\left[\sum_{i=1}^m\hspace{-0.1cm} \|x_i^{k+1}-\bar x^{k+1}\|^2|\mathcal{F}^k\hspace{-0.05cm}\right]
\hspace{-0.1cm}\le \hspace{-0.1cm} (1+a^k -\kappa\gamma^k)\hspace{-0.1cm}\sum_{i=1}^m\hspace{-0.1cm} \|x_i^k - \bar x^k\|^2 + \b^k
\end{equation}
where
$\b^k=a^k\left( \|\bar x^k - \theta^\ast\|^2+\sum_{i=1}^m \|x_i^k - \bar x^k\|^2\right)$.
 Since
$\sum_{k=0}^\infty a^k<\infty$ \as by our assumption, and
the sequences $\{\|\bar x^k-\theta^*\|^2\}$
and $\{\sum_{i=1}^m \|x_i^k-\bar x^k\|^2\}$ are bounded {\it a.s.}, it follows that $\sum_{k=0}^\infty\b^k<\infty$ \as
Thus, the preceding relation satisfies the conditions of
Lemma~\ref{Lemma-polyak_2} with
$v^k= \sum_{i=1}^m \|x_i^k - \bar x^k\|^2$, $q^k=\kappa\gamma^k$  and $p^k=\b^k$
due to our assumptions $\sum_{k=0}^\infty b^k<\infty$ \as and
$\sum_{k=0}^\infty \g^k=\infty$.
So one yields \as
\be\label{eq-sumable}
 \sum_{k=0}^\infty \kappa\gamma^k\sum_{i=1}^m \|x_i^k - \bar x^k\|^2<\infty,\:
 \lim_{k\to\infty} \sum_{i=1}^m \|x_i^k - \bar x^k\|^2=0
 \ee

 It remains to show that $\|\bar x^k-\theta^*\|^2\to0$ \as
For this,
we consider relation~\eqref{eq:Theorem_decreasing}
and focus on the first element of $\bv^k$,  for which we obtain \as for all $k\ge0$:
\begin{equation}\label{eq:x_opt}
\begin{aligned}
&\mathbb{E}\left[\|\bar x^{k+1}-\theta^*\|^2|\mathcal{F}^k\right]
\le  (1+a^k)\|\bar x^k -\theta^*\|^2 \cr
&+\left(\frac{\gamma^k}{m}+a^k\right) \sum_{i=1}^m \|x_i^k - \bar x^k\|^2 + b^k -c^k  (F(\bar x^k)- F(\theta^*))
\end{aligned}
\end{equation}
The preceding relation satisfies Lemma~\ref{lem-opt}
with $\phi=F$, $z^*=\theta^*$,
$z^k =\bar x^k$, $\a^k=a^k$, $\eta^k= c^k$, and
$\beta^k=
(\frac{\gamma^k}{m}+a^k)
\sum_{i=1}^m \|x_i^k - \bar x^k\|^2+b^k$.
By our assumptions, the sequences $\{a^k\}$ and $\{b^k\}$ are summable {\it a.s.}, and $\sum_{k=0}^\infty c^k=\infty$. In view of~\eqref{eq-sumable}, it follows that
$\sum_{k=0}^\infty\beta^k<\infty$ \as
Hence, all the conditions of Lemma~\ref{lem-opt}
are satisfied and, consequently, $\{\bar x^k\}$ converges \as to some optimal solution.

\subsection{Proof of Theorem \ref{theorem:convergence_algorithm_1}}
The basic idea is to apply
Lemma \ref{th-main_decreasing} to the quantities
$\mathbb{E}\left[\|\bar x^{k+1}-\theta^*\|^2|\mathcal{F}^k\right]$ and $\mathbb{E}\left[\sum_{i=1}^m\|x_i^{k+1}-\bar x^{k+1}\|^2|\mathcal{F}^k\right]$.
We divide the proof into two parts to analyze $\|\bar x^{k+1}-\theta^*\|^2$ and $\sum_{i=1}^m\|x_i^{k+1}-\bar x^{k+1}\|^2$, respectively.

Part I: We first analyze  $\|\bar x^{k+1}-\theta^*\|^2$.
For the sake of notational simplicity, we represent $\nabla f_i(x_i^k)$ as $g_i^k$. Stacking  $x_i^k$  and $g_i^k$ into augmented vectors $(x^k)^T=\left[(x_1^k)^T,\cdots,(x_m^k)^T\right]$ and $(g^k)^T=\left[(g_1^k)^T,\cdots,(g_m^k)^T\right]$, respectively, we can write the dynamics of
Algorithm 1~as
\begin{equation}\label{eq:augmented_x_Algorithm}
x^{k+1}=(I+\gamma^k W\otimes I_d)x^k+\gamma^k \zeta_w^k-\lambda^k g^k
\end{equation}
where $\otimes$ denotes the Kronecker product, and $( \zeta_w^k)^T=\left[(\zeta_{w1}^k)^T,\cdots,(\zeta_{wm}^k)^T\right]$ with
$\zeta_{wi}^k\triangleq \sum_{j\in\mathbb{N}_i^{\rm in}}w_{ij}\zeta_j^k$.

From~\eqref{eq:augmented_x_Algorithm} we can obtain the following relationship for the average vector $\bar x^k=\frac{1}{m}\sum_{i=1}^m x_i^k$
\begin{equation}\label{eq:bar_x^k}
\bar x^{k+1}=\bar x^k +\gamma^k \bar\zeta_w^k- \frac{\lambda^k }{m}\sum_{i=1}^m g_i^k
\end{equation}
where
$\bar{\zeta}_w^k=\frac{1}{m}\sum_{i=1}^{m}\zeta_{wi}^k=\frac{1}{m}\sum_{i=1}^{m}\sum_{j\in\mathbb{N}_i^{\rm in}}w_{ij}\zeta_j^k =-\frac{\sum_{i=1}^mw_{ii}\zeta_i^k}{m}$ (note $w_{ii}\triangleq-\sum_{j\in\mathbb{N}^{\rm in}_i}w_{ij}$).

Using (\ref{eq:bar_x^k}) and the preceding relation, we relate $\bar{x}^k$ to an optimal solution
\begin{equation*}
\bar x^{k+1}-\theta^{\ast}=\bar x^k-\theta^{\ast}-\frac{1}{m}\sum_{i=1}^m\left(\lambda^k g_i^k +\gamma^kw_{ii}\zeta_i^k\right)
\end{equation*}
which further implies
\begin{equation*}
\begin{aligned}
&\left\|\bar x^{k+1}-\theta^*\right\|^2 \hspace{-0.1cm}=\hspace{-0.1cm}\left\|\bar
x^k -\theta^* \right\|^2\hspace{-0.1cm} -\hspace{-0.1cm} \frac{2}{m} \hspace{-0.1cm}\sum_{i=1}^m
\hspace{-0.1cm}\left\langle\hspace{-0.05cm} \lambda^k g_i^k+\gamma^kw_{ii}\zeta_i^k , \bar x^k -\theta^*
\hspace{-0.05cm}\right\rangle  \\
 & \qquad\qquad+ \frac{1}{m^2} \left\|\sum_{i=1}^m\left(\lambda^k g_i^k+\gamma^kw_{ii}\zeta_i^k\right)\right\|^2\\
 &\qquad\leq \left\|\bar x^k -\theta^* \right\|^2 - \frac{2}{m} \sum_{i=1}^m
\left\langle \lambda^k g_i^k+\gamma^kw_{ii}\zeta_i^k , \bar x^k -\theta^*
\right\rangle  \\
 &\qquad\qquad+  \frac{2}{m^2}\left\|\sum_{i=1}^m \lambda^k g_i^k\right\|^2
 + \frac{2}{m^2}\left\|\sum_{i=1}^m \gamma^kw_{ii}\zeta_i^k\right\|^2
\end{aligned}
\end{equation*}
 Taking the conditional expectation, given $\mathcal{F}^k=\{x^0,\ldots,x^k\}$, and
using the assumption that the noise $\zeta_i^k$ is with zero mean and variance
$(\sigma_{i}^k)^2$ conditionally on $x_i^k$
(see Assumption~\ref{assumption:dp-noise}),
from the preceding relation we obtain \as for all $k\ge0$,
\begin{equation}\label{eq:conditional_expectation_diminishing}
\begin{aligned}
\mathbb{E}&\left[ \left\|\bar x^{k+1}-\theta^*\right\|^2|\mathcal{F}^k\right]\\
&  \leq   \left\|\bar x^k -\theta^* \right\|^2  - \frac{2\lambda^k}{m} \sum_{i=1}^m
\left\langle  g_i^k, \bar x^k -\theta^*
\right\rangle  \\
  & \qquad+  \frac{2}{m^2}(\lambda^k)^2\left\|\sum_{i=1}^m  g_i^k\right\|^2
  + \frac{2}{m}(\gamma^k)^2\sum_{i=1}^m w_{ii}^2(\sigma_{i}^k)^2
\end{aligned}
\end{equation}

We next estimate the inner product term, for which we have
\begin{equation}\label{eq:inner_product_0}
\begin{aligned}
 &\frac{2\lambda^k}{m}\sum_{i=1}^m\left\langle g_i^k, \bar x^k -\theta^* \right\rangle =  \frac{2\lambda^k}{m}\sum_{i=1}^m\left\langle g_i^k-\nabla f_i(\bar{x}^k),  \bar x^k -\theta^*
 \right\rangle \\
 &\qquad+ \frac{2\lambda^k}{m}\sum_{i=1}^m\left\langle\nabla f_i(\bar{x}^k), \bar x^k -\theta^*
 \right\rangle
 \end{aligned}
\end{equation}
Recalling that $g_i^k=\nabla f_i(x_i^k)$,
by the Lipschitz continuous property of $\nabla f_i(\cdot)$, we have
\begin{equation}\label{eq:inner_product_1}
\begin{aligned}
&\lambda^k\left\langle g_i^k-\nabla f_i(\bar{x}^k),  \bar x^k-\theta^*\right\rangle \geq -L\lambda^k\|x_i^k-\bar{x}^k\|\|\bar x^k -\theta^*\|\\
 &\geq- \frac{\gamma^k}{2}\|x_i^k-\bar{x}^k\|^2
- \frac{ L^2 (\lambda^k)^2 }{2\gamma^k}\|\bar x^k -\theta^*\|^2
\end{aligned}
\end{equation}
By the convexity of $F(\cdot)$, we have
\begin{equation}\label{eq:inner_product_2}
\begin{aligned}
\frac{2\lambda^k}{m}\sum_{i=1}^m\left\langle\nabla f_i(\bar{x}^k), \bar x^k -\theta^*
 \right\rangle&= 2\lambda^k \left\langle\nabla F(\bar{x}^k), \bar x^k -\theta^*
 \right\rangle\\
 &\geq 2\lambda^k (F(\bar{x}^k)-F(\theta^*))
 \end{aligned}
\end{equation}
Combining   (\ref{eq:inner_product_0}), (\ref{eq:inner_product_1}), and (\ref{eq:inner_product_2}) leads to
\begin{equation}\label{eq:inner_product_final}
\begin{aligned}
& \frac{2\lambda^k}{m}\sum_{i=1}^m\left\langle g_i^k, \bar x^k -\theta^* \right\rangle
\geq -\frac{\gamma^k}{m}\sum_{i=1}^m\|x_i^k-\bar{x}^k\|^2
\\
&\qquad-  \frac{L^2(\lambda^k)^2}{\gamma^k}\|\bar x^k -\theta^*\|^2
+2\lambda^k (F(\bar{x}^k)-F(\theta^{\ast}))
 \end{aligned}
\end{equation}
We next estimate the  second last term in (\ref{eq:conditional_expectation_diminishing}):
\begin{equation}\label{eq:second_last_item}
\begin{aligned}
\frac{1}{m^2} \left\|\sum_{i=1}^m g_i^k\right\|^2
&=  \frac{1}{m^2} \left\|\sum_{i=1}^m \left(g_i^k-\nabla f_i(\theta^{\ast})\right)\right\|^2\\
&\leq \frac{L^2}{m} \sum_{i=1}^m\left\| x_i^k- \theta^{\ast} \right\|^2=\frac{L^2}{m} \| x^k- x^{\ast}\|^2
\end{aligned}
\end{equation}

Further using the inequality
\begin{equation}\label{eq:x^k-barx_norm}
\begin{aligned}
\|x^k-x^{\ast}\|^2&\leq \|x^k- {\bf 1}\otimes\bar x^k+ {\bf
1}\otimes\bar x^k-x^{\ast}\|^2\\
&\leq 2\|x^k- {\bf 1}\otimes\bar x^k\|^2 +2\|{\bf 1}\otimes\bar
x^k-x^{\ast}\|^2\\
&\leq 2\sum_{i=1}^m \|x_i^k-\bar{x}^k\|^2 +2m\| \bar
x^k-\theta^{\ast}\|^2
\end{aligned}
\end{equation}
we have from (\ref{eq:second_last_item}) that
\begin{equation}\label{eq:second_last_item_final}
\begin{aligned}
\frac{1}{m^2} \left\|\sum_{i=1}^m g_i^k\right\|^2
\leq \frac{2L^2}{m} \sum_{i=1}^m \|x_i^k-\bar{x}^k\|^2
+ 2L^2  \| \bar x^k-\theta^{\ast}\|^2
\end{aligned}
\end{equation}

Substituting (\ref{eq:inner_product_final}) and (\ref{eq:second_last_item_final})
into  (\ref{eq:conditional_expectation_diminishing})  yields
\begin{equation}\label{eq:bar_x_and_theta}
\begin{aligned}
&\mathbb{E}\left[ \left\|\bar x^{k+1}-\theta^*\right\|^2|\mathcal{F}^k\right]
 \leq   \left\|\bar x^k -\theta^* \right\|^2
 +\frac{\gamma^k}{m}\sum_{i=1}^m\|x_i^k-\bar{x}^k\|^2
\\
&\quad+  L^2(\lambda^k)^2\left(\frac{1}{\gamma^k}+4\right)
\|\bar x^k -\theta^*\|^2
-2\lambda^k (F(\bar{x}^k)-F(\theta^{\ast}))\\
&\quad+\frac{4L^2(\lambda^k)^2}{m} \sum_{i=1}^m \|x_i^k-\bar{x}^k\|^2
+  \frac{2(\gamma^k)^2}{m}\sum_{i=1}^m w_{ii}^2(\sigma_{i}^k)^2
\end{aligned}
\end{equation}

Part II: Next we  analyze $\sum_{i=1}^m\|x_i^{k+1}-\bar x^{k+1}\|^2$.
Using (\ref{eq:augmented_x_Algorithm}) and (\ref{eq:bar_x^k}), we obtain
\[
\begin{aligned}
&x^{k+1}-{\bf 1}\otimes \bar{x}^{k+1}=(I+\gamma^k W\otimes I_d)x^k-{\bf 1}\otimes \bar{x}^{k}\\
&\quad+\gamma^k\left(\zeta_w^k-\frac{1}{m}\sum_{i=1}^m {\bf 1}\otimes \zeta^k_{w,i} \right)-\lambda^k \left(g^k  -\frac{1}{m}\sum_{i=1}^m {\bf 1}\otimes g_i^k\right)
\end{aligned}
\]
Noting  ${\bf 1}\otimes \bar{x}^{k}=\frac{1}{m}\left({\bf 1}{\bf 1}^T\otimes I_d\right) x^k$, $\sum_{i=1}^m{\bf 1}\otimes   \zeta_{w,i}^{k}=\left({\bf 1}{\bf 1}^T\otimes I_d\right) \zeta_w^k$, and $\sum_{i=1}^m{\bf 1}\otimes g_i^{k}=\left({\bf 1}{\bf 1}^T\otimes I_d\right)g^k$,
we can rewrite the preceding equality as
\begin{equation}\label{eq:x-barx}
\begin{aligned}
&x^{k+1}-{\bf 1}\otimes \bar{x}^{k+1} = \hat{W}_k x^k  +\gamma^k\Xi \zeta_w^k -\lambda^k \Xi g^k
\end{aligned}
\end{equation}
with $\hat{W}_k\triangleq \left(I+\gamma^k W- \frac{1}{m}{\bf 1}{\bf 1}^T\right)\otimes I_d$
and $\Xi\triangleq \left(I-\frac{1}{m} {\bf 1}{\bf 1}^T \right)\otimes I_d$.

Since   $\left(I+\gamma^kW-\frac{1}{m} {\bf 1}{\bf 1}^T \right){\bf 1}=0$ holds and we always have $(A\otimes B)(C\otimes D)=(AC)\otimes (BD)$, it follows that
\[
\begin{aligned}
&\hat{W}_k\left( {\bf 1}\otimes \bar{x}^k\right) =\left(\left(I+\gamma^kW-\frac{{\bf 1}{\bf 1}^T}{m} \right){\bf \times 1}\right)\otimes\left(I_d\times \bar{x}^k\right)= 0
\end{aligned}
\]
By subtracting  $\hat{W}_k\left( {\bf 1}\otimes \bar{x}^k\right)=0$ from the right hand side of (\ref{eq:x-barx}), we obtain
\begin{equation*}
\begin{aligned}
&x^{k+1}-{\bf 1}\otimes \bar{x}^{k+1}=\hat{W}_k \left(x^k- {\bf 1}\otimes \bar{x}^{k}\right) +\gamma^k\Xi\zeta_w^k -\lambda^k \Xi g^k
\end{aligned}
\end{equation*}
which further leads to
\[\begin{aligned}
&\|x^{k+1}-{\bf 1}\otimes \bar{x}^{k+1}\|^2\\
&=\|\hat{W}_k(x^k- {\bf 1}\otimes \bar{x}^{k})-\lambda^k \Xi g^k\|^2 + \|\gamma^k \Xi \zeta_w^k\|^2 \\
&\quad+ 2\left\langle \hat{W}_k(x^k- {\bf 1}\otimes \bar{x}^{k})-\lambda^k \Xi g^k,\gamma^k \Xi \zeta_w^k\right \rangle\\
&\leq  \|\hat{W}_k(x^k- {\bf 1}\otimes \bar{x}^{k})-\lambda^k \Xi g^k\|^2 + m(\gamma^k)^2  \sum_{i=1}^m\sum_{j\in\mathbb{N}_i^{\rm in}} w_{ij}^2\|\zeta^k_j \|^2\\
&\quad + 2\left\langle \hat{W}_k(x^k- {\bf 1}\otimes \bar{x}^{k})-\lambda^k \Xi g^k,\gamma^k \Xi  \zeta_w^k\right \rangle
\end{aligned}
\]
where the  inequality follows from $\|\Xi\|=1$ and the definition $\zeta_{wi}^k\triangleq \sum_{j\in\mathbb{N}_i^{\rm in}}w_{ij}\zeta_j^k$.
Taking the conditional expectation with respect to
 $\mathcal{F}^k=\{x^0,\ldots,x^k\}$ and using Assumption~\ref{assumption:dp-noise}  yield
\[
\begin{aligned}
&\mathbb{E}\left[  \|x^{k+1}-{\bf 1}\otimes \bar{x}^{k+1}\|^2|\mathcal{F}^k\right]\\
&\leq \left\|\hat{W}_k(x^k- {\bf 1}\otimes \bar{x}^{k})-\lambda^k \Xi g^k\right\|^2\hspace{-0.15cm}
+ m(\gamma^k)^2\hspace{-0.1cm}   \sum_{i=1}^m\sum_{j\in\mathbb{N}_i^{\rm in}}w_{ij}^2(\sigma_{j}^k)^2 \\
&\hspace{-0.1cm}\leq\left(\|\hat{W}_k(x^k\hspace{-0.05cm}- \hspace{-0.05cm}{\bf 1}\otimes \bar{x}^{k})\|\hspace{-0.05cm}+\hspace{-0.05cm}\|\lambda^k \Xi g^k\|\right)^2\hspace{-0.15cm} +\hspace{-0.05cm}m(\gamma^k)^2\max_{j\in[m]} (\sigma^k_{j} )^2 C_W
\end{aligned}
\]
where
$C_W=\sum_{i=1}^m\sum_{j\in\mathbb{N}_i^{\rm in}}w_{ij}^2$.
Using the fact $\|\Xi\|=1$ and $\|\hat{W}_k\|=\|I+\gamma^kW-\frac{1}{m}{\bf  1}{\bf  1}^T\|=1-\gamma^k|\nu|$ where $-\nu$ is some non-zero eigenvalue of $W$  (see Assumption~\ref{assumption:W}), we obtain
\begin{equation}\label{eq:x^k-barx_norm2}
\hspace{-0.1cm}\begin{aligned}
&\mathbb{E}\left[ \|x^{k+1}-{\bf 1}\otimes \bar{x}^{k+1}\|^2|\mathcal{F}^k\right]\\
&  \leq (1-\gamma^k|\nu|)^2(1+\epsilon) \|x^{k}-{\bf 1}\otimes \bar{x}^{k}\|^2  \cr
&\quad+ (1+\epsilon^{-1})(\lambda^k)^2\|g^k\|^2
+ m(\gamma^k)^2\max_{j\in[m]}(\sigma^k_{j})^2 C_W
\end{aligned}
\end{equation}
for any $\epsilon>0$, where  we used
$(a+b)^2\le (1+\epsilon) a^2 + (1+\epsilon^{-1})b^2$ valid for any scalars $a,b,$ and $\epsilon>0$.

We next focus on estimating the term involving the gradient $g^k$ in the preceding inequality.
Noting   $g^k=m\nabla f(x^k)$  and that $f(\cdot)$ has Lipschitz continuous gradients (with Lipschitz constant $\frac{L}{m}$), we have
\[
\begin{aligned}
\begin{aligned}
\|g^k\|^2 &=m^2\|\nabla f(x^k)-\nabla f(x^*)+\nabla f(x^*)\|^2 \\
&\leq 2m^2\|\nabla f(x^k)-\nabla f(x^*)\|^2+2m^2\|\nabla f(x^*)\|^2\\
&\le  2L^2\|x^k- x^*\|^2+2m^2\|\nabla f(x^*)\|^2
\end{aligned}
\end{aligned}
\]
Since $x^*={\bf 1}\otimes\theta^*$, using the relationship in (\ref{eq:x^k-barx_norm}), we obtain
\[
\|g^k\|^2 \leq 4L^2(\|x^k- {\bf 1}\otimes \bar{x}^k\|^2+m\|\bar{x}^k- \theta^*\|^2)+2m^2\|\nabla f(x^*)\|^2
\]
Finally,  substituting the preceding relation back in (\ref{eq:x^k-barx_norm2}) yields
\[
\begin{aligned}
&\mathbb{E}\left[ \|x^{k+1}-{\bf 1}\otimes \bar{x}^{k+1}\|^2|\mathcal{F}^k\right]
\\&\leq (1-\gamma^k|\nu|)^2(1+\epsilon) \|x^{k}-{\bf 1}\otimes \bar{x}^{k}\|^2\\
&+ 4(1+\epsilon^{-1})L^2(\lambda^k)^2(\|x^k- {\bf 1}\otimes \bar{x}^k\|^2+m\|\bar{x}^k- \theta^*\|^2)  \\
&+ 2(1+\epsilon^{-1})(\lambda^k)^2 m^2\|\nabla f(x^*)\|^2 +m(\gamma^k)^2
\max_{j\in[m]}(\sigma^k_{j} )^2 C_W
\end{aligned}
\]
By letting $\epsilon=\frac{\gamma^k|\nu|}{1- \gamma^k|\nu|}$ 
 and consequently $1+\epsilon=(1-\gamma^k|\nu|)^{-1}$ and $1+\epsilon^{-1}=(\gamma^k|\nu|)^{-1}$,
we arrive at
\begin{equation}\label{eq:x^k-barx_norm_final}
\begin{aligned}
&\mathbb{E}\left[ \|x^{k+1}-{\bf 1}\otimes \bar{x}^{k+1}\|^2|\mathcal{F}^k\right]\\
&\leq 
\left(1-\gamma^k|\nu| +\frac{4L^2(\lambda^k)^2}{|\nu|\gamma^k}\right) \|x^{k}-{\bf 1}\otimes \bar{x}^{k}\|^2\\
& +  \frac{4 m L^2(\lambda^k)^2}{|\nu|\gamma^k} \|\bar{x}^k- \theta^*\|^2 +\frac{4(\lambda^k)^2 m^2}{|\nu|\gamma^k}\|\nabla f(x^*)\|^2 \\
&  + m(\gamma^k)^2\max_{j\in[m]}(\sigma^k_{j} )^2 C_W
\end{aligned}
\end{equation}
By combining~(\ref{eq:bar_x_and_theta}) and (\ref{eq:x^k-barx_norm_final}), and using   Assumption~\ref{assumption:dp-noise}, we  have $\mathbb{E}\left[\|\bar x^{k+1}-\theta^*\|^2|\mathcal{F}^k\right]$ and $\mathbb{E}\left[\sum_{i=1}^m\|x_i^{k+1}-\bar x^{k+1}\|^2|\mathcal{F}^k\right]$ satisfying the conditions of Lemma~\ref{th-main_decreasing} with
$\kappa=|\nu|$, $c^k=2\lambda^k$,
$a^k=\max\{L^2(\lambda^k)^2\left(\frac{1}{\gamma^k}+4\right),\frac{4 m L^2(\lambda^k)^2}{|\nu|\gamma^k}\}$, and $b^k=(\gamma^k)^2\max\{ \frac{2}{m}\sum_{i=1}^m w_{ii}^2(\sigma_{i}^k)^2, \, \frac{4(\lambda^k)^2 m^2}{|\nu|\gamma^k}\|\nabla f(x^*)\|^2+m \max_{j\in[m]}(\sigma^k_{j} )^2 C_W\}$ where
$C_W =\sum_{i=1}^m\sum_{j\in\mathbb{N}_i^{\rm in}}w_{ij}^2$.


\subsection{Proof of Lemma \ref{Theorem:general_gradient_tracking}}

Since the results of Lemma~\ref{th-dsystem} are asymptotic, they remain valid when the starting index is shifted from $k=0$ to $k=T$, for an arbitrary $T\ge 0$. So the  idea is to show that the conditions in Lemma~\ref{th-dsystem}
are satisfied  for all $k\ge T$ (for some large enough $T\geq 0$).

\noindent(a)  Because $\kappa_i>0$ for all $1\leq i\leq 9$, for $\pi=[\pi_1, \pi_2, \pi_3]^T$ to satisfy  $\pi^T V\leq \pi^T$ and $\pi^TH^k\geq 0$,  we only need to show that the following inequalities can be true
\begin{equation}\label{eq:conditons_pi}
\begin{aligned}
 &\kappa_1\lambda^k\pi_1+ (1-\kappa_3\gamma_1^k) \pi_2 \leq \pi_2,\\
& \kappa_2\lambda^k\pi_1+(1-\kappa_4 \gamma_2^k) \pi_3\leq \pi_3,
\\
 &\left(\kappa_6 \lambda^k-\kappa_7(\lambda^k)^2\right)\pi_1- \kappa_8\frac{(\lambda^k)^2}{\gamma_1^k}\pi_2-\kappa_9\frac{(\lambda^k)^2}{\gamma_2^k}\pi_3\geq 0
\end{aligned}
\end{equation}

The first inequality is equivalent to $\pi_2\geq \frac{\kappa_1\lambda^k}{\kappa_3\gamma_1^k}\pi_1$. Given that $\lim_{k\to\infty}\lambda^k/\gamma_1^k=0$ holds and $\gamma_1^k$ as well as $\lambda^k$ are positive   according to the assumption, it can easily be seen that for any given $\pi_1>0$, we can always find a $\pi_2>0$ satisfying the relationship when $k$ is larger than some $T \geq 0$.

The second inequality is equivalent to $\pi_3\geq \frac{\kappa_2\lambda^k}{\kappa_4\gamma_2^k}\pi_1$. Given that $\lim_{k\to\infty}\lambda^k/\gamma_2^k=0$ holds and $\gamma_2^k$ as well as $\lambda^k$ are positive   according to the assumption, it can easily be seen that for any given $\pi_1>0$, we can always find a $\pi_3>0$ satisfying the relationship when $k$ is larger than some $T \geq 0$.

The third inequality is equivalent to $\pi_1\geq \frac{\kappa_7\lambda^k}{\kappa_6}\pi_1+\frac{\kappa_8\lambda^k}{\kappa_6\gamma_1^k}\pi_2+\frac{\kappa_9\lambda^k}{\kappa_6\gamma_2^k}\pi_3$. Since the right hand side converges to zero according to our assumptions on $\lambda^k$, $\gamma_1^k$ and $\gamma_2^k$, we can always find a constant $\pi_1$ satisfying this inequality for $k\geq T$.
 Thus,   we can always find a vector $\pi$  satisfying all  inequalities in (\ref{eq:conditons_pi}) for $k \ge T$ for some large enough $T\geq 0$, and hence the conditions in Lemma~\ref{th-dsystem} are satisfied.

 By Lemma~\ref{th-dsystem}, it follows that
for the three entries of $\bv^k$, i.e., $\bv_1^k$, $\bv_2^k$, and
$\bv_3^k$, we have that
\be\label{eq-limit-exist} \lim_{k\to\infty}
\pi_1\bv^k_1+\pi_2 \bv_2^k+\pi_3\bv_3^k \ee
 exists \as and
$\sum_{k=0}^\infty\pi^T H^k\bu^k<\infty$ holds \as with $\bu^k= [\|\nabla
F(\bar x^k)\|^2,\ \|\bar y^k\|^2]^T$. Since
$
 \pi^TH^k= \left[\kappa_5\lambda^k\pi_1,\
(\kappa_6 \lambda^k-\kappa_7(\lambda^k)^2)\pi_1-\kappa_8\frac{(\lambda^k)^2}{\gamma_1^k}\pi_2-\kappa_9\frac{(\lambda^k)^2}{\gamma_2^k}\pi_3\right]$
and  $(\lambda^k)^2$, $\frac{(\lambda^k)^2}{\gamma_1^k}$, and $\frac{(\lambda^k)^2}{\gamma_2^k}$ are summable,
 one has \be\label{eq-sumfinite}
\sum_{k=0}^\infty\lambda^k\|\nabla F(\bar x^k)\|^2<\infty, \quad
\sum_{k=0}^\infty\lambda^k\|\bar y^k\|^2<\infty, \quad a. s.\ee
 Hence, it follows that
\be\label{eq-to0}
 \| \nabla F(\bar x^k)\|<\Delta_1,\quad
 \|\bar y^k\|<\Delta_2\quad a.s.
\ee
for some random scalars $\Delta_1>0$ and $\Delta_2>0$ due to the assumption
$\sum_{k=0}^\infty\lambda^k=\infty$.

Now, we focus on proving that both $\bv_2^k= \|{\bf x}^{k}-{\bf \bar x}^{k}\|_R^2$ and $\bv_3^k=
\|{\bf y}^{k}-{\rm diag}(v){\bf \bar y}^{k}\|_C^2$ converge \as to
0. The  idea is to show that we can apply Lemma~\ref{Theo:convergence_to_zero}.
By focusing on the second and third elements of $\bv^k$, i.e., $\bv_2^k$ and $\bv_3^k$, from
\eqref{eq-fine} we have
\bas \left[\begin{array}{c} \bv_2^{k+1}\cr
\bv_3^{k+1}\end{array} \right] \le \left(\tilde V^k + a^k {\bf 1}{\bf
1}^T\right) \left[\begin{array}{c} \bv_2^k\cr
\bv_3^k\end{array}\right] +\hat b^k{\bf 1} + \left[\begin{array}{c}
\hat{c}^k \cr
\hat{c}^k\end{array}\right]
\eas
where
$\hat b^k = b^k+a^k(F(\bar x^k)-F^*)$, $
 \hat{c}^k =\max\left\{\kappa_8\frac{(\lambda^k)^2}{\gamma_1^k}\|\bar y^k\|^2,\, \kappa_9\frac{(\lambda^k)^2}{\gamma_2^k}\|\bar y^k\|^2\right\}
$, and $\tilde
V^k=\left[\begin{array}{cc}
 1-\kappa_3\gamma_1^k & 0\cr
0 & 1-\kappa_4\gamma_2^k\end{array}\right]$, which can be rewritten as
\ba\label{eq-finer}
\left[\begin{array}{c} \bv_2^{k+1}\cr
\bv_3^{k+1}\end{array} \right] \le \tilde V^k \left[\begin{array}{c}
\bv_2^k\cr \bv_3^k\end{array}\right] +\tilde b^k{\bf 1} \ea
 where
$\tilde b^k= b^k+\hat{c}^k+a^k\left( F(\bar x^k) - F^*  + \|{\bf x}^{k}-{\bf \bar x}^{k}\|_R^2+
\|{\bf y}^{k}-{\rm diag}(v){\bf \bar y}^{k}\|_C^2\right)$.

To apply Lemma \ref{Theo:convergence_to_zero}, noting that $\g_1^k$ and $\g_2^k$ are not summable, we show that the equation
 $\tilde \pi^T
\tilde V^k=(1-\a\gamma_1^k) \tilde\pi^T$ has a solution in $\tilde\pi=[\pi_2,\pi_3]$  with
$\pi_2,\pi_3>0$ and $\a\in(0,1)$.
From $\tilde \pi^T
\tilde V^k=(1-\a\gamma_1^k) \tilde \pi^T$, one has
\[
(1-\kappa_3\gamma_1^k)\pi_2 \leq (1-\alpha\gamma_1^k)\pi_2,\quad (1-\kappa_4\gamma_2^k)\pi_3\leq (1-\alpha\gamma_1^k)\pi_3
\]
which can be simplified as
$
\alpha\leq \kappa_3$, $\alpha\leq \frac{\kappa_4\g_2^k}{\g_1^k}
$.

Given   $\lim_{k\to\infty}\gamma_1^k/\gamma_2^k<\infty$ according to our assumption, it can be seen that $\frac{\kappa_4\g_2^k}{\g_1^k}$ is positive, and hence, such an $\alpha\in (0,1)$ and $\tilde\pi>0$ can always be found.

 We next  prove that the condition $\sum_{k=0}^{\infty}\tilde{b}^k<0$ \as
of Lemma \ref{Theo:convergence_to_zero} is
also   satisfied. Indeed, the condition can be met because: (1) $b^k$, $a^k$, $\frac{(\lambda^k)^2}{\gamma_1^k}$, and $\frac{(\lambda^k)^2}{\gamma_2^k}$ are all summable according to the assumption of the lemma; and (2)   $\|\bar y^k\|$ (see (\ref{eq-to0}))  and $F(\bar x^k) -
F^*$, $\|{\bf x}^{k}-{\bf \bar x}^{k}\|_R^2$, $\|{\bf y}^{k}-{\rm diag}(v){\bf \bar y}^{k}\|_C^2$ are all bounded \as
due to  the existence of the limit in (\ref{eq-limit-exist}).
Thus, all the conditions of
Lemma~\ref{Theo:convergence_to_zero} are satisfied, so
it follows that  $\lim_{k\to\infty}\|x_i^k - \bar
x^k\|=0$ and $\lim_{k\to\infty}\|y_i^k - v_i\bar y^k\|=0$ \as
 Moreover, in view of the existence of the limit in (\ref{eq-limit-exist}) and the facts that $\pi_1>0$ and $v_1^k=F(\bar x^k)-F(\theta^*)$, it follows that $\lim_{k\to\infty} F(\bar x^k)$
exists \as

\noindent(b) Since $\sum_{k=0}^\infty\lambda^k\|\nabla F(\bar
x^k)\|^2<\infty$ holds \as (see~\eqref{eq-sumfinite}), from $\sum_{k=0}^\infty
\lambda^k =\infty$, it follows that we have
$\lim\inf_{k\to\infty}\|\nabla F(\bar x^k)\|=0$ \as

 Now, if the function $F(\cdot)$ has bounded level sets, then the sequence
$\{\bar x^k\}$ is \as bounded since   $\lim_{k\to\infty} F(\bar x^k)$
exists \as (as shown in part (a)). Thus, $\{\bar x^k\}$ \as
has accumulation points.
Let $\{\bar x^{k_i}\}$ be a sub-sequence such that
$\lim_{i\to\infty}\|\nabla F(\bar x^{k_i})\|=0$ \as
Without loss of
generality, we may assume that $\{\bar x^{k_i}\}$ is \as convergent, for
otherwise we would choose a sub-sequence of $\{\bar x^{k_i}\}$. Let
$\lim_{i\to\infty}\bar x^{k_i}=\hat x$. Then, by the continuity of the
gradient $\nabla F(\cdot)$, it follows $\nabla F(\hat x)=0$,
implying that $\hat x$ is an optimal point. Since $F(\cdot)$ is continuous,
 we have  $\lim_{i\to\infty} F(\bar x^{k_i}) = F(\hat x)=F^*$.
By part (a),  $\lim_{k\to\infty} F(\bar x^k)$
exists {\it a.s.}, so we must have $\lim_{k\to\infty} F(\bar x^k)=F^*$ \as

Finally,  by part (a), we  have $\lim_{k\to\infty}\|x_i^k-\bar
x^k\|^2=0$ \as for every $i$. Thus,  each  $\{x_i^k\}$
has the same accumulation points as the sequence $\{\bar x^k\}$ {\it a.s.}, implying by
the continuity of the function $F(\cdot)$ that
 $\lim_{k\to\infty}F(x^k_i)=F^*$ \as for all $i$.

\subsection{Proof of Theorem \ref{theorem_convergence_algorithm2}}
We divide the derivations into four
steps: in Step I, Step II, and Step III, we establish  relations
for $\mathbb{E}\left[ F(\bar x^k)-F^*|\mathcal{F}^k\right]$,
$\mathbb{E}\left[\|{\bf x}^{k}-{\bf \bar x}^{k}\|_R^2|\mathcal{F}^k\right]$, and
$\mathbb{E}\left[\|{\bf y}^{k}-{\rm diag}(v){\bf \bar y}^{k}\|_C^2|\mathcal{F}^k\right]$ for the iterates generated by Algorithm 2, respectively. In Step IV, we use them to
show that~(\ref{eq-fine}) of Lemma~\ref{Theorem:general_gradient_tracking} holds.

Step I: Relationship for $\mathbb{E}\left[F(\bar x^{k}) -F^*|\mathcal{F}^k\right]$.

Since $F$ is convex with Lipschitz gradients, we have
\[F(y)\le F(x)+\la \nabla F(x),y-x\ra+\frac{L}{2}\|y-x\|^2, \: \forall y,x\in\mathbb{R}^d\]
Letting $y=\bar x^{k+1}$ and $x=\bar x^k$ in the preceding relation
and using  (\ref{eq:x_bar_update2}), as well as
$\nabla F(\bar{x}^k)=\frac{1}{m}\sum_{i=1}^m \nabla f_i(\bar x^k)$, we obtain
\[
\begin{aligned}
&F(\bar x^{k+1}\hspace{-0.07cm})
\le F(\bar x^k)\hspace{-0.07cm}-\hspace{-0.07cm} \left\la\hspace{-0.1cm} \nabla F(\bar{x}^k),\lambda^k \frac{(u\otimes {  I_d})^T}{m}\hspace{-0.07cm}\left(y^k\hspace{-0.07cm}-(v\otimes {I_d})\bar{y}^k\right)\right.\\
&\left.+\lambda^k\frac{u^Tv}{m}\bar{y}^k-\gamma_1^k \bar{\zeta}_w^k
\right\ra\\
&+\frac{L}{2}\left\|\lambda^k \frac{(u\otimes {  I_d})^T}{m}\left(y^k-(v\otimes {I_d})\bar{y}^k\right)+\lambda^k\frac{u^Tv}{m}\bar{y}^k-\gamma_1^k \bar{\zeta}_w^k
\right\|^2
\end{aligned}
\]
We estimate the last term in the preceding inequality by using  $(a+b+c)^2\le 3a^2 + 3b^2+3c^2$ and $\|A\otimes B\|_2=\|A\|_2\|B\|_2$.
After subtracting $F^*$ on both sides of the resulting inequality, we obtain
\[
\begin{aligned}
&F(\bar x^{k+1})-F^{\ast}\le  F(\bar x^k)-F^{\ast}\\
&\hspace{-0.07cm}- \left\la \hspace{-0.1cm}\nabla F(\bar{x}^k),\lambda^k \frac{(u\otimes {  I_d})^T}{m}\hspace{-0.1cm}\left(y^k\hspace{-0.1cm}-\hspace{-0.1cm}(v\otimes {I_d})\bar{y}^k\right) \hspace{-0.1cm}+\hspace{-0.1cm}\lambda^k\frac{u^Tv}{m}\bar{y}^k\hspace{-0.1cm}-\hspace{-0.1cm}\gamma_1^k \bar{\zeta}_w^k
\right\ra \\
& + \frac{3L(\lambda^k)^2\|u\|^2 }{2m^2} \left\|y^k-(v\otimes {I_d})\bar{y}^k \right\|^2+ \frac{3L(\lambda^k)^2(u^Tv)^2}{2m^2} \left\|\bar{y}^k \right\|^2  \\
&  +\frac{3L(\gamma_1^k)^2}{2} \left\|\bar\zeta_w^k\right\|^2
\end{aligned}
\]

Taking the conditional expectation  with respect to $\mathcal{F}^k$ yields
\begin{equation}\label{eq:F-F^ast}
\begin{aligned}
&\mathbb{E}\left[ F(\bar x^{k+1})-F^{\ast}|\mathcal{F}^k\right]
\le F(\bar x^k)-F^{\ast}\\
&- \left\la \nabla F(\bar{x}^k),\lambda^k \frac{(u\otimes {  I_d})^T}{m}\left(y^k-(v\otimes {I_d})\bar{y}^k\right)+\lambda^k\frac{u^Tv}{m}\bar{y}^k
\right\ra  \\
&+ \frac{3L(\lambda^k)^2\|u\|^2 }{2m^2} \hspace{-0.07cm}\left\|y^k\hspace{-0.07cm}-(v\otimes {I_d})\bar{y}^k \right\|^2 \hspace{-0.07cm} + \frac{3L(\lambda^k)^2(u^Tv)^2}{2m^2}\hspace{-0.07cm} \left\|\bar{y}^k \right\|^2\\
& + \frac{3L(\gamma_1^k)^2}{2}(\sigma^k_{R,\zeta})^2
\end{aligned}
\end{equation}
where
$
(\sigma^k_{R,\zeta})^2\triangleq  \mathbb{E}\left[\left\|\bar\zeta_w^k\right\|^2\right]
= \frac{1}{m^2} \sum_{j=1}^{m} \left(\textstyle\sum_{i\in\mathbb{N}_{R,j}^{\rm out}}u_i R_{ij}\right)^2(\sigma^k_{\zeta,j})^2
$.

The inner product term  in~\eqref{eq:F-F^ast} satisfies
\begin{equation}\label{eq-inner_product}
\begin{aligned}
&-\left\la \nabla F(\bar{x}^k),\lambda^k \frac{(u\otimes {  I_d})^T}{m}\left(y^k-(v\otimes {I_d})\bar{y}^k\right)+\lambda^k\frac{u^Tv}{m}\bar{y}^k
\right\ra  \\
&=-\lambda^k\frac{u^Tv}{m}\left\la \nabla F(\bar{x}^k), \frac{(u\otimes {  I_d})^T}{u^Tv}\left(y^k-(v\otimes {I_d})\bar{y}^k\right)+\bar{y}^k
\right\ra  \\
& =\frac{u^Tv\lambda^k}{2m}
\left\|\nabla F(\bar{x}^k)-\bar{y}^k  -\frac{(u\otimes {  I_d})^T}{u^Tv}\left(y^k-(v\otimes {I_d})\bar{y}^k\right)\right\|^2
\\
& -\hspace{-0.1cm}\frac{u^T\hspace{-0.1cm}v\lambda^k}{2m}\hspace{-0.1cm}\left\| \nabla F(\bar x^k)
\right\|^2 \hspace{-0.15cm} -\hspace{-0.1cm}\frac{u^T\hspace{-0.1cm}v\lambda^k}{2m}\hspace{-0.1cm}\left\|\frac{\hspace{-0.1cm}(u\hspace{-0.06cm}\otimes {  \hspace{-0.06cm}I_d})^T}{u^Tv}\hspace{-0.1cm}\left(y^k\hspace{-0.1cm}-\hspace{-0.1cm}(v\hspace{-0.06cm}\otimes\hspace{-0.06cm} {I_d})\bar{y}^k\right)\hspace{-0.06cm}+\hspace{-0.06cm}\bar{y}^k\right\|^2
\end{aligned}
\end{equation}
where in the second equality we use $-\la a,b\ra=\frac{\|a-b\|^2-\|a\|^2-\|b\|^2}{2}$ valid for any vectors $a$ and $b$.

Using the relationship
\[
\nabla F(\bar{x}^k)-\bar{y}^k =\frac{1}{m}\sum_{i=1}^m(\nabla f_i(\bar x^k)-g_i^k)+\bar{g}^k-\bar{y}^k
\]
and the inequality $\|a+b+c\|^2\leq 3\|a\|^2+3\|b\|^2+3\|c\|^2$, we can bound the first term on the right hand side of (\ref{eq-inner_product}) as follows:
\[
\begin{aligned}
  & \left\|\nabla F(\bar{x}^k)-\bar{y}^k
  -\frac{(u\otimes {  I_d})^T}{u^Tv}\hspace{-0.1cm}\left(y^k\hspace{-0.1cm}-(v\otimes {I_d})\bar{y}^k\right)\right\|^2  \\
      &\leq 3\left\|\sum_{i=1}^m\frac{\nabla f_i(\bar x^k)-g_i^k}{m}\right\|^2 +3\|\bar{g}^k-\bar{y}^k\|^2\\
      &\qquad+3\left\|\frac{(u\otimes {  I_d})^T}{u^Tv}\left(y^k-(v\otimes {I_d})\bar{y}^k\right)\right\|^2
\end{aligned}
\]

Using the inequality $\|a+b\|^2\geq  \|a\|^2+ \|b\|^2$, we can bound the last term on the right hand side of (\ref{eq-inner_product}) as follows:
\[
\begin{aligned}
&-\frac{u^Tv\lambda^k}{2m}\left\|\frac{(u\otimes {  I_d})^T}{u^Tv}\left(y^k-(v\otimes {I_d})\bar{y}^k\right)+\bar{y}^k\right\|^2\\
&\leq\hspace{-0.07cm} -\frac{u^Tv\lambda^k}{2m}\hspace{-0.07cm}\left\|\frac{(u\otimes {  I_d})^T}{u^Tv}\left(y^k\hspace{-0.1cm}-(v\otimes {I_d})\bar{y}^k\right)\right\|^2\hspace{-0.15cm}-\frac{u^Tv\lambda^k}{2m}\hspace{-0.07cm}\left\|\bar{y}^k\right\|^2
\end{aligned}
\]
Plugging the preceding two inequalities into~(\ref{eq-inner_product}) and combining
the common terms lead to
\begin{equation}\label{eq-inner_product2}
\begin{aligned}
&-\left\la \nabla F(\bar{x}^k),\lambda^k \frac{(u\otimes {  I_d})^T}{m}\left(y^k-(v\otimes {I_d})\bar{y}^k\right)+\lambda^k\frac{u^Tv}{m}\bar{g}^k
\right\ra  \\
&\leq \frac{3u^Tv\lambda^k}{2m}\left\|\sum_{i=1}^m\frac{\nabla f_i(\bar x^k)-g_i^k}{m}\right\|^2+\frac{3u^Tv\lambda^k}{2m}\|\bar{g}^k-\bar{y}^k\|^2 \\
&\qquad+\frac{u^Tv\lambda^k}{m}\left\|\frac{(u\otimes {  I_d})^T}{u^Tv}\left(y^k-(v\otimes {I_d})\bar{y}^k\right)\right\|^2
 \\
&\qquad -\frac{u^Tv\lambda^k}{2m}\left\| \nabla F(\bar x^k)
\right\|^2 -\frac{u^Tv\lambda^k}{2m}\left\|\bar{y}^k\right\|^2\\
&\leq \frac{3u^TvL^2\lambda^k}{2m^2}\sum_{i=1}^m\left\| \bar{x}^k-x_i^k\right\|^2 + \frac{3u^Tv\lambda^k}{2m}\|\bar{g}^k-\bar{y}^k\|^2\\
&\qquad+\frac{\|u\|^2\lambda^k}{ m u^Tv}\left\|  y^k-(v\otimes {I_d})\bar{y}^k \right\|^2
 \\
&\qquad-\frac{u^Tv\lambda^k}{2m}\left\| \nabla F(\bar x^k)
\right\|^2 -\frac{u^Tv\lambda^k}{2m}\left\|\bar{y}^k\right\|^2
%
%
%
\end{aligned}
\end{equation}
where in the last inequality we used $\|A\otimes B\|_2=\|A\|_2\|B\|_2$ and  the assumption that each $\nabla f_i(\cdot)$ has Lipschitz continuous gradients with the constant $L$.

Next, combining (\ref{eq:F-F^ast}) and (\ref{eq-inner_product2}) yields
\begin{equation}\label{eq:F-F^ast_final1}
\begin{aligned}
\begin{aligned}
&\mathbb{E}\left[ F(\bar x^{k+1})-F^{\ast}|\mathcal{F}^k\right]\\
&\le F(\bar x^k)-F^{\ast}+\frac{3u^TvL^2\lambda^k}{2m^2}\sum_{i=1}^m\left\| x_i^k-\bar{x}^k\right\|^2 \\
&+\left(\frac{\|u\|^2\lambda^k}{ m u^Tv} + \frac{3L(\lambda^k)^2\|u \|^2 }{2m^2} \right)\sum_{i=1}^{m}\left\|  y_i^k-v_i\bar{y}^k \right\|^2\\
&
 +\frac{3u^Tv\lambda^k}{2m}\|\bar{g}^k-\bar{y}^k\|^2 -\frac{u^Tv\lambda^k}{2m}\left\| \nabla F(\bar x^k)
\right\|^2  \\
&-\left(\frac{u^Tv\lambda^k}{2m}- \frac{3L(\lambda^k)^2(u^Tv)^2}{2m^2}\right)\left\|\bar{y}^k\right\|^2 +\frac{3L(\gamma_1^k)^2}{2}(\sigma^k_{R,\zeta})^2
\end{aligned}\end{aligned}
\end{equation}

Using the fact that in finite dimensional vector spaces, all norms are equivalent up to a proportionality constant, we  always have a constant $\delta_{2,R}$  such that $\|x\|_2\leq \delta_{2,R}\|x\|_R$ for all $x$ and hence $\sum_{i=1}^m\left\| x_i^k-\bar{x}^k\right\|^2\leq \delta_{2,R}^2 \|{\bf x}^{k}-{\bf \bar x}^{k}\|_R^2$ according to the definition of matrix norm in (\ref{eq:matrix_norm_R}). Similarly, we always have $\sum_{i=1}^{m}\left\|  y_i^k-v_i\bar{y}^k \right\|^2\leq \delta_{2,C}^2  \|{\bf y}^{k}-{\rm diag}(v){\bf \bar y}^{k}\|_C^2$ for some constant $\delta_{2,C}$ according to the definition of matrix norm in (\ref{eq:matrix_norm_C}). Therefore, (\ref{eq:F-F^ast_final1}) can be rewritten as
\begin{equation}\label{eq:F-F^ast_final}
\begin{aligned}
\begin{aligned}
&\mathbb{E}\left[ F(\bar x^{k+1})-F^{\ast}|\mathcal{F}^k\right]\\
& \le F(\bar x^k)-F^{\ast}+\frac{3u^TvL^2\lambda^k\delta_{2,R}^2 \|{\bf x}^{k}-{\bf \bar x}^{k}\|_R^2}{2m^2}  \\
&+\left(\frac{\|u\|^2\lambda^k}{m u^Tv} + \frac{3L(\lambda^k)^2\|u \|^2 }{2m^2} \right)\delta_{2,C}^2  \|{\bf y}^{k}-{\rm diag}(v){\bf \bar y}^{k}\|_C^2\\
&
 +\frac{3u^Tv\lambda^k}{2m}\|\bar{g}^k-\bar{y}^k\|^2 -\frac{u^Tv\lambda^k}{2m}\left\| \nabla F(\bar x^k)
\right\|^2  \\
&-\left(\frac{u^Tv\lambda^k}{2m}- \frac{3L(\lambda^k)^2(u^Tv)^2}{2m^2}\right)\left\|\bar{y}^k\right\|^2  +\frac{3L(\gamma_1^k)^2}{2}(\sigma^k_{R,\zeta})^2
\end{aligned}\end{aligned}
\end{equation}

Step II: Relationship for $\|{\bf x}^{k+1}-{\bf \bar x}^{k+1}\|_R^2$.\\
For the convenience of analysis, we write the iterates of  (\ref{eq:push-pull}) on per-coordinate expressions. Define for all
$\ell=1,\ldots,d,$ and $k\ge0$,
$x^k(\ell)=([x_1^k]_\ell,\ldots,[x_m^k]_\ell)^T$,
$y^k(\ell)=([y_1^k]_\ell,\ldots,[y_m^k]_\ell)^T$, $g^k(\ell)=([g_1^k]_\ell,\ldots,[g_m^k]_\ell)^T$, $\zeta_w^k(\ell) =([\zeta_{w1}^k]_\ell,\ldots,[\zeta_{wm}^k]_\ell)^T$, $\xi_w^k(\ell)=([\xi_{w1}^{k}]_\ell,\ldots,[\xi_{wm}^{k}]_\ell)^T
$. In
this per-coordinate view, (\ref{eq:push-pull})  has the following form  for all $\ell=1,\ldots,d,$
and $k\ge0$,
\begin{equation}\label{eq-alg2-percord}
\begin{aligned}
x^{k+1}(\ell)&=(I+\gamma_1^k R)x^k(\ell)+\gamma_1^k \zeta_w^k(\ell)-\lambda^ky^k(\ell)\cr
y^{k+1}(\ell)&=(I-\alpha^k +\gamma_2^kC)y^k(\ell) +  \gamma_2^k \xi_w^k(\ell)
\\
&\quad+ g^{k+1}(\ell) -(1-\alpha^k)g^k(\ell)
\end{aligned}
\end{equation}

 From the definition of $x^{k+1}(\ell)$ in~\eqref{eq-alg2-percord},
and the relation for  $\bar x^{k+1}$ in
(\ref{eq:bar_x_k+1}), we obtain for all $\ell=1,\ldots,d$,
\[
\begin{aligned}
x^{k+1}(\ell)& - [\bar x^{k+1}]_\ell{\bf 1}
=(I+\gamma_1^k R)\left(x^k(\ell) -[\bar x^k]_\ell{\bf 1}\right)\\
&\quad+\gamma_1^k \left(\zeta_w^k(\ell)- [\bar\zeta_w^k]_\ell{\bf 1}\right) -\lambda^k\left(I -\frac{{\bf 1}u^T}{m}\right)y^k(\ell)
\end{aligned}
\]
where we used $(I+\gamma_1^k R){\bf 1}={\bf 1}$.

Noting that
$[\bar x^k]_\ell$ is the average of $x^k(\ell)$, i.e.,
$\frac{1}{m}{\bf 1}u^T\left(x^k(\ell) -[\bar x^k]_\ell{\bf
1}\right)=0$, and abbreviating $I -\frac{{\bf 1}u^T}{m}$ as $\Pi_u\triangleq I -\frac{{\bf 1}u^T}{m}$, we have
\[
\begin{aligned}
x^{k+1}(\ell)\hspace{-0.06cm} - \hspace{-0.06cm}[\bar x^{k+1}]_\ell{\bf 1}&=\hspace{-0.06cm}\bar R^k\hspace{-0.1cm}\left(x^k(\ell)
\hspace{-0.06cm}-\hspace{-0.06cm}[\bar x^k]_\ell{\bf 1}\right)\hspace{-0.06cm}+\hspace{-0.06cm}\gamma_1^k \left(\zeta_w^k(\ell)\hspace{-0.06cm}-\hspace{-0.06cm} [\bar\zeta_w^k]_\ell{\bf 1}\right) \\
&\quad-\lambda^k\Pi_u \left(y^k(\ell)-v[\bar{y}^k]_\ell\right) -\lambda^k\Pi_uv[\bar{y}^k]_\ell
\end{aligned}
\]
where $\bar{R}^k=I+\gamma_1^kR-\frac{ {\bf 1}u^T}{m} $.

Taking norm $\|\cdot\|_R$ on both sides   leads to
\[
\begin{aligned}
&\left\|x^{k+1}(\ell) - [\bar x^{k+1}]_\ell{\bf 1} \right\|_R^2=\\
&\left\|\bar R^k\hspace{-0.1cm}\left(x^k(\ell)
-[\bar x^k]_\ell{\bf 1}\right) \hspace{-0.1cm}-\hspace{-0.1cm}\lambda^k\Pi_u \left(y^k(\ell)-v[\bar{y}^k]_\ell\right)\hspace{-0.1cm}-\hspace{-0.1cm}\lambda^k\Pi_uv[\bar y^k]_\ell\right\|_R^2 \\
&+\left\| \gamma_1^k \left(\zeta_w^k(\ell)- [\bar\zeta_w^k]_\ell{\bf 1}\right)  \right\|_R^2 \\
&+ 2\left\langle R^k\left(x^k(\ell)
-[\bar x^k]_\ell{\bf 1}\right) -\lambda^k\Pi_u \left(y^k(\ell)-v[\bar{y}^k]_\ell\right)\right.\\
&\left. \qquad-\lambda^k\Pi_uv[\bar y^k]_\ell,\gamma_1^k \left(\zeta_w^k(\ell)- [\bar\zeta_w^k]_\ell{\bf 1}\right) \right\rangle_R
\end{aligned}
\]
 where $\la\cdot\ra_R$ denotes the inner product induced\footnote{It can be seen that the norm $\|\cdot\|_R$ satisfies the Parallelogram Law and, hence, the inner product induced by the norm $\|\cdot\|_R$ exists.} by the norm $\|\cdot\|_R$.

Taking the expectation (conditioned on $\mathcal{F}^k$)   yields
\begin{equation}\label{eq:x_k+1-bar_x_k+1_element}
\begin{aligned}
&\mathbb{E}\left[\left\|x^{k+1}(\ell) - [\bar x^{k+1}]_\ell{\bf 1} \right\|_R^2|\mathcal{F}^k\right]\\
&=\left\|\bar R^k\hspace{-0.1cm}\left(x^k(\ell)\hspace{-0.1cm}
-\hspace{-0.1cm}[\bar x^k]_\ell{\bf 1}\right)\hspace{-0.1cm} -\hspace{-0.1cm}\lambda^k\Pi_u \left(y^k(\ell)\hspace{-0.1cm}-\hspace{-0.1cm}v[\bar{y}^k]_\ell\right)\hspace{-0.1cm}-\hspace{-0.1cm}\lambda^k\Pi_u v[\bar y^k]_\ell\right\|_R^2\\
&\qquad +\mathbb{E}\left[\left\| \gamma_1^k \left(\zeta_w^k(\ell)- [\bar\zeta_w^k]_\ell{\bf 1}\right)  \right\|_R^2\right]\\
&\leq \left(\left(1-\gamma_1^k\rho_R\right)\left\|  x^k(\ell)
-[\bar x^k]_\ell{\bf 1} \right\|_R \right.\\
&\qquad \left.+\left\|\lambda^k\Pi_u \left(y^k(\ell)-v[\bar{y}^k]_\ell\right)+\lambda^k\Pi_uv[\bar y^k]_\ell\right\|_R\right)^2\\
&\qquad +\mathbb{E}\left[\left\| \gamma_1^k \left(\zeta_w^k(\ell)- [\bar\zeta_w^k]_\ell{\bf 1}\right)  \right\|_R^2\right]\\
&\leq \hspace{-0.1cm}\left(1\hspace{-0.06cm}-\hspace{-0.06cm}\gamma_1^k\rho_R\right)\left\|  x^k(\ell)
\hspace{-0.06cm}-\hspace{-0.06cm}[\bar x^k]_\ell{\bf 1} \right\|_R^2  +\frac{2(\lambda^k)^2\left\|\Pi_uv\right\|_R^2}{\rho_R\gamma_1^k}\left\| [\bar y^k]_\ell\right\|_R^2\\
 & +\frac{2(\lambda^k)^2\left\|\Pi_u\right\|_R^2}{\rho_R\gamma_1^k}\left\| y^k(\ell)-v[\bar{y}^k]_\ell\right\|^2_R\\
 &  +(\gamma_1^k)^2 \mathbb{E}\left[\left\| \zeta_w^k(\ell)- [\bar\zeta_w^k]_\ell{\bf 1} \right\|_R^2\right]
\end{aligned}
\end{equation}
where in the last equality we used the relationship
$(a+b)^2\le (1+\epsilon) a^2 + (1+\epsilon^{-1})b^2$  on the first term on the right hand side of (\ref{eq:x_k+1-bar_x_k+1_element}) by setting $\epsilon$ to $\frac{1}{1-\gamma_1^k\rho_R}-1$ (resulting in $1+\frac{1}{\epsilon}=\frac{1}{\rho_R\gamma_1^k}$).
Summing the
preceding relations over $\ell=1,\ldots,d$, and noting
$\sum_{\ell=1}^d\|x^{k+1}(\ell)-[\bar x^{k+1}]_\ell{\bf 1}\|_R^2
=\|{\bf x}^{k+1}-{\bf \bar x}^{k+1}\|_R^2$,
,  $\sum_{\ell=1}^d\|y^k(\ell)-v[\bar y^k]_\ell \|_R^2=\|{\bf y}^{k}-{\rm diag}(v){\bf \bar y}^{k}\|_R^2$, and $\sum_{\ell=1}^d\left\| \zeta_w^k(\ell)- [\bar\zeta_w^k]_\ell{\bf 1} \right\|_R^2\leq \delta_{R,2}^2 \sum_{i=1}^m\left\| \zeta_w^k-  \bar\zeta_w^k  \right\|_2^2$, we obtain
\begin{equation}\label{eq:x_k+1-bar_x_k+1_final}
\begin{aligned}
&\mathbb{E}\left[\|{\bf x}^{k+1}-{\bf \bar x}^{k+1}\|_R^2|\mathcal{F}^k\right]\\
&\leq   \left(1-\gamma_1^k\rho_R\right)\|{\bf x}^{k}-{\bf \bar x}^{k}\|_R^2 +\frac{2(\lambda^k)^2\left\|\Pi_uv\right\|_R^2}{\rho_R\gamma_1^k}\left\|\bar y^k \right\|_R^2 \\
&\quad+\frac{2(\lambda^k)^2\left\|\Pi_u\right\|_R^2}{\rho_R\gamma_1^k}\|{\bf y}^{k}-{\rm diag}(v){\bf \bar y}^{k}\|_R^2\\
&\quad +2m(\gamma_1^k)^2\delta_{R,2}^2(\sigma_{R,\zeta}^k)^2 + 2(\gamma_1^k)^2\delta_{R,2}^2\sum_{i,j}(R_{ij}\sigma_{\zeta,j}^k)^2
\end{aligned}
\end{equation}
where $\delta_{R,2}$ is constant such that $\|x\|_R\leq \delta_{R,2}\|x\|_2$ for all $x$. (In finite dimensional vector spaces, all norms are equivalent up to a proportionality constant,  represented by $\delta_{R,2}$ here.)

Step III: Relationship for $\|{\bf y}^{k}-{\rm diag}(v){\bf \bar y}^{k}\|_C^2$.\\
From
(\ref{eq:bar_x_k+1}), the $\ell$th entries of  $[\bar y^k]_\ell$
satisfy
\[[\bar y^{k+1}]_\ell
=(1-\alpha^k)[\bar y^k]_\ell + \gamma_2^k[\bar{\xi}_w^{k}]_\ell+ [\bar{g}^{k+1}]_\ell-
 (1-\alpha^k)[\bar{g}^k]_\ell\]
Then, using
 \eqref{eq-alg2-percord}, we obtain for all $\ell=1,\ldots,d$,
\[
\begin{aligned}
y^{k+1}(\ell) - v[\bar y^{k+1}]_\ell   =& \bar C^k(y^k(\ell)-v[\bar y^k]_\ell ) +\gamma_2^k\Pi_v \xi_w^{k}(\ell)
\\
&+\Pi_v \left(g^{k+1}(\ell)- (1-\alpha^k)g^{k}(\ell)\right)
\end{aligned}
\]
where $\bar C^k=(1-\alpha^k)(I-\frac{1}{m} {v}{\bf 1}^T)+\gamma_2^kC$ and $
\Pi_v=\left(I-\frac{1}{m}v{\bf 1}^T\right)$. Note that we used the relationship $ \bar{C}^k  v[\bar y^k]_\ell=0$.

Taking the norm $\|\cdot\|_C$ on both sides yields
\[
\begin{aligned}
&\left\|y^{k+1}(\ell) - v[\bar y^{k+1}]_\ell \right\|_C^2 \\
=&\left\| \bar C^k(y^k(\ell)-v[\bar y^k]_\ell)+\Pi_v \left(g^{k+1}(\ell) -(1-\alpha^k) g^k(\ell)\right)\right\|_C^2 \\
&+ (\gamma_2^k)^2\left\|\Pi_v \left(\xi_w^{k}(\ell)\right) \right\|_C^2+ 2\left\la \bar C^k(y^k(\ell)-v[\bar y^k]_\ell)\right.\\
&\left.+\Pi_v \left(g^{k+1}(\ell) -(1-\alpha^k) g^k(\ell)\right)  , \gamma_2^k\Pi_v  \xi_w^{k}(\ell)\right\ra_C\\
\leq& \left((1-\rho_C\gamma_2^k)\left\| y^k(\ell)-v[\bar y^k]_\ell \right\|_C\right.\\
&\left.+\|\Pi_v\|_C\left\|g^{k+1}(\ell) - (1-\alpha^k)g^k(\ell)\right\|_C\right)^2 \\
&+ (\gamma_2^k)^2\|\Pi_v\|_C^2\left\|\xi_w^{k}(\ell) \right\|_C^2\\
&+ 2\left\la \bar C^k(y^k(\ell)-v[\bar y^k]_\ell )+\Pi_v \left(g^{k+1}(\ell) - (1-\alpha^k)g^k(\ell)\right),\right.\\
&\left. \gamma_2^k\Pi_v  \xi_w^{k}(\ell)\right\ra_C
\end{aligned}
\]
where in the inequality we used $\|\bar{C}^k\|_C\leq 1-\gamma_2^k\rho_C$.
Taking conditional expectation on both sides leads to
\begin{equation}\label{eq:y^k+1-bar_y^k+1}
\begin{aligned}
&\mathbb{E}\left[\left\|y^{k+1}(\ell) - v[\bar y^{k+1}]_\ell \right\|_C^2|\mathcal{F}^k\right]\\
&\leq  (1-\gamma_2^k\rho_C) \left\|  y^k(\ell)-v[\bar y^k]_\ell\right\|_C^2\\
&\qquad +\frac{\|\Pi_v\|_C^2}{ \gamma_2^k\rho_C }\mathbb{E}\left[\left\| g^{k+1}(\ell) - (1-\alpha^k)g^k(\ell) \right\|_C^2|\mathcal{F}^k\right] \\
&\qquad+ (\gamma_2^k)^2\|\Pi_v\|_C^2\mathbb{E}\left[\left\| \xi_w^{k}(\ell) \right\|_C^2\right]
\end{aligned}
\end{equation}
where   we used the relationship  $(a+b)^2\le (1+\epsilon) a^2 + (1+\epsilon^{-1})b^2$ valid for any scalars $a,b$ and $\epsilon>0$ and set $\epsilon= \frac{1}{1-\gamma_2^k\rho_C}-1$.

By summing these relations over $\ell=1,\ldots,d$, we find
\begin{equation}\label{eq:y^k+1-bar_y^k+1_2}
\begin{aligned}
\mathbb{E}&\left[\|{\bf y}^{k+1}-{\rm diag}(v){\bf \bar y}^{k+1}\|_C^2|\mathcal{F}^k\right]\\
&\leq  (1-\gamma_2^k\rho_C) \|{\bf y}^{k}-{\rm diag}(v){\bf \bar y}^{k}\|_C^2 \\
&\qquad+\frac{\|\Pi_v\|_C^2\delta^2_{C,2}}{\gamma_2^k\rho_C}\mathbb{E}\left[\sum_{i=1}^m\left\| g_i^{k+1} - (1-\alpha^k)g_i^k \right\|_2^2|\mathcal{F}^k\right] \\
&\qquad+(\gamma_2^k)^2\|\Pi_v\|_C^2\delta^2_{C,2}\mathbb{E}\left[\sum_{i=1}^m\left\|\xi_{wi}^{k}\right\|_2^2\right]
\end{aligned}
\end{equation}

Using Assumption \ref{as:bounded_gradients}, we have
\[
\begin{aligned}
&\left\| g_i^{k+1} - (1-\alpha^k)g_i^k\right\|\\
&\leq \left\| g_i^{k+1} - g_i^k\right\|+\|\alpha^kg_i^k\|\leq L\left\| x_i^{k+1} - x_i^k\right\|+\alpha^kC,
\end{aligned}
\]
which, in combination with (\ref{eq:y^k+1-bar_y^k+1_2}), yields
\begin{equation}\label{eq:y^k+1-bar_y^k+1_3}
\begin{aligned}
&\mathbb{E}\left[\|{\bf y}^{k+1}-{\rm diag}(v){\bf \bar y}^{k+1}\|_C^2|\mathcal{F}^k\right]\\
&\leq  (1-\gamma_2^k\rho_C) \|{\bf y}^{k}-{\rm diag}(v){\bf \bar y}^{k}\|_C^2\\
&\quad+\frac{2L^2\|\Pi_v\|_C^2\delta^2_{C,2}}{\gamma_2^k\rho_C}\mathbb{E}\left[\sum_{i=1}^m\left\| x_i^{k+1} - x_i^k \right\|_2^2|\mathcal{F}^k\right] \\
&\quad+ \frac{2(\alpha^k)^2\|\Pi_v\|_C^2\delta^2_{C,2}C^2}{\gamma_2^k\rho_C}+ (\gamma_2^k)^2 \|\Pi_v\|_C^2\delta^2_{C,2}\sum_{i,j} (C_{ij}\sigma_{\xi,j}^k)^2
\end{aligned}
\end{equation}
where  $\delta_{C,2}$ is constant such that $\|x\|_C\leq \delta_{C,2}\|x\|_2$ for all $x$. (In finite dimensional vector spaces, all norms are equivalent up to a proportionality constant, represented by $\delta_{C,2}$ here.)

Next, we proceed to analyze $\mathbb{E}\left[\sum_{i=1}^m\left\| x_i^{k+1} - x_i^k\right\|_2^2|\mathcal{F}^k\right]$ in (\ref{eq:y^k+1-bar_y^k+1_2}). Using (\ref{eq:push-pull}), we have for every  index
$\ell=1,\ldots,d$:
\[
\begin{aligned}
& x^{k+1}(\ell)-x^k(\ell) =\gamma_1^kRx^k(\ell)+\gamma_1^k\zeta_w^k(\ell)-\lambda^ky^k(\ell)\\
&  \hspace{-0.06cm}=\hspace{-0.06cm}\gamma_1^kR(x^k\hspace{-0.04cm}(\ell)\hspace{-0.06cm}-\hspace{-0.06cm}[\bar x^k]_\ell{\bf
1})\hspace{-0.06cm}+\hspace{-0.06cm}\gamma_1^k\hspace{-0.04cm}\zeta_w^k\hspace{-0.04cm}(\ell)\hspace{-0.06cm}-\hspace{-0.06cm}\lambda^k\hspace{-0.06cm}(y^k\hspace{-0.04cm}(\ell)\hspace{-0.06cm}-\hspace{-0.06cm}v[\bar{y}^k]_\ell  ) \hspace{-0.06cm}-\hspace{-0.06cm}\lambda^k \hspace{-0.04cm}v[\bar{y}^k]_\ell
\end{aligned}
\] where we used
  $R{\bf 1}=0$ in the second equality.

The preceding relationship leads to
\[
\begin{aligned}
&\mathbb{E}\left[\left\|x^{k+1}(\ell)-x^k(\ell) \right\|_2^2 |\mathcal{F}^k\right]\\
&\leq   \left(r^k\left\|x^k(\ell)-[\bar x^k]_\ell{\bf
1}\right\|_2+\lambda^k\left\|y^k(\ell)-v[\bar{y}^k]_\ell \right\|_2+\right.\\
&\quad \left.\lambda^k\|v\|_2\left\|[\bar{y}^k]_\ell \right\|_2\right)^2  + \mathbb{E}\left[\left\|\gamma_1^k\zeta_w^k(\ell) \right\|_2^2\right]\\
&\leq 3(r^k)^2\left\|x^k(\ell)-[\bar x^k]_\ell{\bf
1}\right\|_2^2+3(\lambda^k)^2\left\|y^k(\ell)-v[\bar{y}^k]_\ell]\right\|_2^2  \\
&\qquad+3 (\lambda^k)^2 \|v\|_2^2 \left\|[\bar{y}^k]_\ell \right\|_2^2 + (\gamma_1^k)^2\mathbb{E}\left[\left\|\zeta_w^k(\ell)\right\|_2^2\right]
\end{aligned}
\]
where $r^k=\|\gamma_1^kR\|_2=\gamma_1^k\rho_{c}$ which is arbitrarily close to the spectral radius of the matrix $\gamma_1^kR$ .


By summing over  $\ell=1,\ldots,d$, we obtain
\begin{equation}\label{eq:y_inequality_final}
\begin{aligned}
&\mathbb{E}\sum_{i=1}^m\left[\left\|x_i^{k+1}-x_i^k \right\|_2^2 |\mathcal{F}^k\right]\\
&\leq 3(r^k)^2\sum_{i=1}^m\left\|x_i^k -\bar x^k  \right\|_2^2+3(\lambda^k)^2\sum_{i=1}^m\left\|y_i^k -v_i\bar{y}^k\right\|_2^2\\
&\qquad+ 3(\lambda^k)^2\|v\|_2^2\left\|\bar{y}^k\right\|_2^2+  (\gamma_1^k)^2 \sum_{i,j} R^2_{ij} (\sigma^k_{\zeta,j})^2
\end{aligned}
\end{equation}

Plugging (\ref{eq:y_inequality_final}) into (\ref{eq:y^k+1-bar_y^k+1_3}) and using $r^k= \gamma_1^k\rho_{c}$ lead  to
\begin{equation}\label{eq:y-bar_y_final}
\begin{aligned}
&\mathbb{E}\left[\|{\bf y}^{k+1}-{\rm diag}(v){\bf \bar y}^{k+1}\|_C^2|\mathcal{F}^k\right]\leq  \\
&\hspace{-0.1cm}\left(\hspace{-0.1cm}1-\gamma_2^k\rho_C+\frac{6L^2\|\Pi_v\|_C^2\delta^2_{C,2}\delta^2_{2,C}(\lambda^k)^2}{\rho_C\gamma_2^k}\right)\hspace{-0.1cm}\|{\bf y}^{k}\hspace{-0.1cm}-{\rm diag}(v){\bf \bar y}^{k}\|_C^2\\
&  +
\frac{6 L^2 \|\Pi_v\|_C^2\delta^2_{C,2}\delta^2_{2,R}\rho^2_{c}(\gamma_1^k)^2}{\rho_C\gamma_2^k}\|{\bf x}^{k}-{\bf \bar x}^{k}\|_R^2\\
& +\frac{6\|\Pi_v\|_C^2\|v\|_C^2L^2\delta^2_{C,2}(\lambda^k)^2}{\rho_C\gamma_2^k} \left\|\bar{y}^k\right\|_2^2\\
&+\frac{2(\alpha^k)^2\|\Pi_v\|_C^2\delta^2_{C,2}C^2}{\gamma_2^k\rho_C} + (\gamma_2^k)^2 \|\Pi_v\|_C^2\delta^2_{C,2}\sum_{i,j} (C_{ij}\sigma_{\xi,j}^k)^2 \\
& +  \frac{2 L^2\delta^2_{C,2} \|\Pi_v\|_2^2 (\gamma_1^k)^2 \sum_{i,j} R^2_{ij} (\sigma^k_{\zeta,j})^2 }{\gamma_2^k\rho_C}
\end{aligned}
\end{equation}

Step IV: We combine Steps I-III and prove the theorem.

Defining $\bv^k=\big[(F(\bar
x^{k+1})-F^*),\|{\bf x}^{k}-{\bf \bar x}^{k}\|_R^2,\|{\bf y}^{k}-{\rm diag}(v){\bf \bar y}^{k}\|_C^2\big]^T$, we
have the following relations from (\ref{eq:F-F^ast_final}),
(\ref{eq:x_k+1-bar_x_k+1_final}), and (\ref{eq:y-bar_y_final}):
\begin{equation}\label{eq:stacked}
\begin{aligned}
\mathbb{E}\left[\bv^{k+1}|\mathcal{F}^k\right]\leq(V^k+A^k)\bv^k-H^k\left[\begin{array}{c} \left\|\nabla
F(\bar x^k)\right\|^2\cr \|\bar y^k\|^2
\end{array}\right]+B^k
\end{aligned}
\end{equation}
where
\[
\begin{aligned}
&V^k=\left[\begin{array}{ccc} 1 &
\frac{3(u^Tv)L^2\delta^2_{2,R}\lambda^k}{2m^2}&  \frac{\|u\|^2\delta^2_{2,C}\lambda^k}{m u^Tv}  \cr
0& 1-\gamma_1^k\rho_R&
0\cr 0&0
 &
1-\gamma_2^k\rho_C
\end{array}\right],\\
&A^k\hspace{-0.08cm}=\hspace{-0.08cm}\left[\begin{array}{ccc} 0 &
0&  \frac{3L\delta^2_{2,C}(\lambda^k)^2\|u \|^2 }{2m^2} \cr
 0& 0&
\frac{2\left\|\Pi_u\right\|_R^2\delta^2_{C,R}(\lambda^k)^2}{\rho_R\gamma_1^k}\cr
 0&
\frac{6 L^2 \|\Pi_v\|_C^2\delta^2_{C,2}\delta^2_{2,R}\rho^2_{c}(\gamma_1^k)^2}{\rho_C\gamma_2^k}&
 \frac{6L^2\|\Pi_v\|_C^2\delta^2_{C,2}\delta^2_{2,C}(\lambda^k)^2}{\rho_C\gamma_2^k}
\end{array}\right],\\
&H^k=\left[\begin{array}{cc}   \frac{u^Tv\lambda^k}{2m} &
\frac{u^Tv\lambda^k}{2m}- \frac{3L(\lambda^k)^2(u^Tv)^2}{2m^2}\cr 0&-\frac{2\left\|\Pi_uv\right\|_R^2\delta^2_{2,R}(\lambda^k)^2}{\rho_R\gamma_1^k}\cr
0&-\frac{6\|\Pi_v\|_C^2\|v\|_C^2L^2\delta^2_{C,2} (\lambda^k)^2}{\rho_C\gamma_2^k}
\end{array}\right],\: B^k=\left[\begin{array}{c} b_1^k
\cr
b_2^k
\cr
b_3^k
\end{array}\right]
\end{aligned}
\]
with $b_1^k=\frac{3L(\gamma_1^k)^2}{2}(\sigma^k_{R,\zeta})^2 +\frac{3u^Tv\lambda^k}{2m}\|\bar{g}^k-\bar{y}^k\|^2 $, $b_2^k=2m(\gamma_1^k)^2\delta_{R,2}^2(\sigma_{R,\zeta}^k)^2 + 2(\gamma_1^k)^2\delta_{R,2}^2\sum_{i,j}(R_{ij}\sigma_{\zeta,j}^k)^2$ $b_3^k=\frac{2(\alpha^k)^2\|\Pi_v\|_C^2\delta^2_{C,2}C^2}{\gamma_2^k\rho_C} + (\gamma_2^k)^2 \|\Pi_v\|_C^2\delta^2_{C,2}\sum_{i,j} (C_{ij}\sigma_{\xi,j}^k)^2   +  \frac{2 L^2\delta^2_{C,2} \|\Pi_v\|_2^2 (\gamma_1^k)^2 \sum_{i,j} R^2_{ij} (\sigma^k_{\zeta,j})^2 }{\gamma_2^k\rho_C}$.

From (\ref{eq:bar_x_k+1}), we have
\[
\begin{aligned}
\bar{y}^{k+1}-\bar{g}^{k+1}&=(1-\alpha^k)(\bar{y}^k-\bar{g}^k) +\gamma_2^k \bar{\xi}_w^k
\end{aligned}
\]
which further implies
\[
\begin{aligned}
\|\bar{y}^{k+1}-\bar{g}^{k+1}\|^2=&(1-\alpha^k)^2\| \bar{y}^k-\bar{g}^k\|^2 +(\gamma_2^k)^2 \|\bar{\xi}_w^k\|^2\\
&+2\left\langle (1-\alpha^k)(\bar{y}^k-\bar{g}^k), \gamma_2^k\bar{\xi}_w^k\right\rangle
\end{aligned}
\]
and
\[
\begin{aligned}
\mathbb{E}&\left[\|\bar{y}^{k+1}-\bar{g}^{k+1}\|^2|\mathcal{F}^k\right]\\
&= (1-2\alpha^k+(\alpha^k)^2)\| \bar{y}^k-\bar{g}^k\|^2 +(\gamma_2^k)^2  \mathbb{E}\left[\|\bar{\xi}_w^k\|^2\right]
\end{aligned}
\]
Given that $(\gamma_2^k)^2  \mathbb{E}\left[\|\bar{\xi}_w^k\|^2\right]$ is summable according to the theorem statement, and $\{\alpha^k\}$ is not summable but square summable, we have that $\|\bar{y}^{k+1}-\bar{g}^{k+1}\|^2$ satisfies the condition in Lemma \ref{Lemma-polyak_2}. Therefore,  the sequence $\alpha^k\|\bar{y}^{k+1}-\bar{g}^{k+1}\|^2$ is summable {\it a.s.} according to Lemma \ref{Lemma-polyak_2}, and hence, $\lambda^k\|\bar{y}^{k+1}-\bar{g}^{k+1}\|^2$ is summable {\it a.s.} under the theorem condition $\lim_{k\to\infty}\lambda^k/\alpha^k<\infty$.

Using Assumption \ref{assumption:dp-noises-intrack},  and the conditions that $(\gamma_1^k)^2$, $(\gamma_2^k)^2$, $\frac{(\alpha^k)^2}{\gamma_2^k}$, and $\frac{(\gamma_1^k)^2}{\gamma_2^k}$ are summable in the theorem statement,  it follows  that all entries of
the matrix $B^k$  are \as summable.  By  defining $\hat{b}^k$ as  the maximum element of $B^k$, we have $B^k\leq \hat{b}^k {\bf 1}$. Therefore,
  $\mathbb{E}\left[ F(\bar x^k)-F^*|\mathcal{F}^k\right]$,
$\mathbb{E}\left[\|{\bf x}^{k}-{\bf \bar x}^{k}\|_R^2|\mathcal{F}^k\right]$, and
$\mathbb{E}\left[\|{\bf y}^{k}-{\rm diag}(v){\bf \bar y}^{k}\|_C^2|\mathcal{F}^k\right]$ for the iterates generated by Algorithm 2 satisfy the conditions of Lemma \ref{Theorem:general_gradient_tracking} and,  hence, the results of Lemma \ref{Theorem:general_gradient_tracking} hold.

\bibliographystyle{IEEEtran}

\bibliography{reference1}
\vspace{-1.2cm}

\end{document}